\pdfoutput=1
\documentclass[a4paper,11pt]{amsart}

\def\definetac{\newif\iftac}

\ifx%
  \tactrue\undefined%
  \definetac%
  \ifx%
    \state\undefined\tacfalse%
  \else%
    \tactrue%
  \fi%
\fi%

\iftac\else\usepackage{amsthm}\fi

\usepackage{
   amssymb
  ,amsmath
  ,stmaryrd
  ,mathrsfs
  ,enumitem
  ,xcolor
  ,datetime
  ,leftidx
  ,mathtools
  ,tikz
  ,url
  ,bbding
  ,newclude
  ,xspace
  ,mathbbol}

\usetikzlibrary{calc}
\def\len{0.53033008589}

\definecolor{darkgreen}{rgb}{0,0.45,0}

\usepackage[cmtip,all,2cell,pdf]{xy}\UseAllTwocells

\usepackage[%
   hyperfootnotes=true
  ,colorlinks%
  ,citecolor=darkgreen%
  ,linkcolor=darkgreen%
  ]{hyperref}

\usepackage[color=gray!20]{todonotes}

\usepackage{xcolor}% provides \colorlet
\usepackage{fixme}

\usepackage[%
 margin=1in
,includehead
,includefoot
,headheight=0.5in
,headsep=0.2in
,footskip=3em
,heightrounded
,left=3.5cm
,right=3.5cm
,top=3.5cm
,bottom=6cm%
]{geometry}

\usepackage[utf8]{inputenc}
\usepackage{xparse}

\ExplSyntaxOn
\NewDocumentCommand{\makeabbrev}{mmm}
 {
  \yoruk_makeabbrev:nnn { #1 } { #2 } { #3 }
 }

\cs_new_protected:Npn \yoruk_makeabbrev:nnn #1 #2 #3
 {
  \clist_map_inline:nn { #3 }
   {
    \cs_new_protected:cpn { #2 } { #1 { ##1 } }
   }
 }
\ExplSyntaxOff

\makeabbrev{\mathbf}{b#1}{b,c,d,e,g,h,i,j,k,l,m,n,o,p,q,r,t,u,v,w,x,y,z,%
              B,C,D,E,G,H,I,J,K,L,M,N,O,P,Q,R,T,U,V,W,X,Y,Z}

\makeabbrev{\boldsymbol}{bs#1}{%
    a,b,c,d,e,g,h,i,j,k,l,m,n,o,p,q,r,s,t,u,v,w,x,y,z,%
    A,B,C,D,E,G,H,I,J,K,L,M,N,O,P,Q,R,S,T,U,V,W,X,Y,Z}

\makeabbrev{\mathsf}{sf#1}{a,b,c,d,e,f,g,h,i,j,k,l,m,n,o,p,q,r,s,t,u,v,w,x,y,z,%
                           A,B,C,D,E,F,G,H,I,J,K,L,M,N,O,P,Q,R,S,T,U,V,W,X,Y,Z}

\makeabbrev{\mathfrak}{f#1}{a,b,c,d,e,f,g,h,j,k,l,m,n,o,p,q,r,s,t,u,v,w,x,y,z,%
                             A,B,C,D,E,F,G,H,I,J,K,L,M,N,O,P,Q,R,S,T,U,V,W,X,Y,Z}

\makeabbrev{\mathcal}{c#1}{A,B,C,D,E,F,G,H,I,J,K,L,M,N,O,P,Q,R,S,T,U,V,W,X,Y,Z}
\makeabbrev{\mathbf}{l#1}{A,B,C,D,E,F,G,H,I,J,K,L,M,N,O,P,Q,R,S,T,U,V,W,X,Y,Z}
\makeabbrev{\mathbb}{s#1}{A,B,C,D,E,F,G,H,I,J,K,L,M,N,O,P,Q,R,S,T,U,V,W,X,Y,Z}

\makeatletter
\newif\ifhyperref
\@ifpackageloaded{hyperref}{\hyperreftrue}{\hyperreffalse}
\iftac
  \let\your@state\state
  \def\state#1{\gdef\currthmtype{#1}\your@state{#1}}
  \let\your@staterm\staterm
  \def\staterm#1{\gdef\currthmtype{#1}\your@staterm{#1}}
  \let\defthm\newtheorem
  \def\currthmtype{}
  \ifhyperref
    \def\autoref#1{\ref*{label@name@#1}~\ref{#1}}
  \else
    \def\autoref#1{\ref{label@name@#1}~\ref{#1}}
  \fi
  \AtBeginDocument{%
    \let\old@label\label%
    \def\label#1{%
      {\let\your@currentlabel\@currentlabel%
        \edef\@currentlabel{\currthmtype}%
        \old@label{label@name@#1}}%
      \old@label{#1}}
  }
\else
  \ifhyperref
    \def\defthm#1#2{%
      \newtheorem{#1}{#2}[section]%
      \expandafter\def\csname #1autorefname\endcsname{#2}%
      \expandafter\let\csname c@#1\endcsname\c@thm}
  \else
    \def\defthm#1#2{\newtheorem{#1}[thm]{#2}}
    \ifx\SK@label\undefined\let\SK@label\label\fi
    \let\old@label\label
    \let\your@thm\@thm
    \def\@thm#1#2#3{\gdef\currthmtype{#3}\your@thm{#1}{#2}{#3}}
    \def\currthmtype{}
    \def\label#1{{\let\your@currentlabel\@currentlabel\def\@currentlabel%
        {\currthmtype~\your@currentlabel}%
        \SK@label{#1@}}\old@label{#1}}
    \def\autoref#1{\ref{#1@}}
  \fi
\fi
\makeatother
\numberwithin{equation}{section}
\newtheorem{thm}{Theorem}[section]
\defthm{cor}{Corollary}
\defthm{prop}{Proposition}
\defthm{lem}{Lemma}
\defthm{sch}{Scholium}
\defthm{assume}{Assumption}
\defthm{claim}{Claim}
\defthm{conj}{Conjecture}
\defthm{hyp}{Hypothesis}
\defthm{fact}{Fact}
\defthm{quest}{Question}
\iftac\theoremstyle{plain}\else\theoremstyle{definition}\fi
\defthm{defn}{Definition}
\defthm{defn-prop}{Definition-Proposition}
\defthm{notn}{Notation}
\iftac\theoremstyle{plain}\else\theoremstyle{remark}\fi
\defthm{rmk}{Remark}
\defthm{eg}{Example}
\defthm{egs}{Examples}
\defthm{ex}{Exercise}
\defthm{ceg}{Counterexample}
\defthm{con}{Construction}
\defthm{warn}{Warning}
\defthm{digr}{Digression}
\defthm{per}{Perspective}
\defthm{disc}{Discussion}
\defthm{term}{Terminology}
\defthm{heur}{Heuristics}

\DeclareSymbolFont{bbold}{U}{bbold}{m}{n}
\usepackage{eucal}
  \DeclareSymbolFontAlphabet{\mathbbb}{bbold}
  \newcommand{\bbone}{{[0]}}
  \newcommand{\bbtwo}{{[1]}}
  \newcommand{\bbthree}{{[2]}}

\DeclareMathOperator{\id}{\mathrm{id}}

\DeclareMathOperator{\lan}{Lan}

\DeclareMathOperator{\Mnd}{Mnd}
\DeclareMathOperator*{\holim}{holim}

\let\xto\xrightarrow

\newcommand{\Nearrow}{\rotatebox[origin=c]{45}{$\Rightarrow$}}

\newcommand{\Swarrow}{\rotatebox[origin=c]{225}{$\Rightarrow$}}

\def\toiso{\xto{\smash{\raisebox{-.5mm}{$\scriptstyle\sim$}}}}

\newcommand{\D}{\sD}
\newcommand{\E}{\sE}
\newcommand{\M}{\sM}

\newcommand{\DeclareCategory}[1]{\textsf{\upshape #1}}
  
  \def\Cat{\DeclareCategory{Cat}}
  \def\cCat{\Cat}
  \def\CAT{\DeclareCategory{CAT}}
  \def\cCAT{\CAT}
  \def\PDER{\DeclareCategory{PDer}}
  \def\cPDER{\PDER}
  \def\cDER{\DeclareCategory{Der}}

\def\class#1{\mathcal{#1}}
  \def\A{\class{A}}

\def\Der{\mathsf{Der}}
\def\Dia{\mathsf{Dia}}
\def\Set{\mathsf{Set}}

\def\dia{\mathrm{dia}}
\def\ho{\mathrm{Ho}}

\def\op{^\mathrm{op}}
\newcommand{\Nat}{\mathrm{Nat}}
\newcommand{\pr}{\mathrm{pr}}
\newcommand{\bDelta}{\boldsymbol{\Delta}}

\def\cSp{\mathrm{Sp}}

\def\Tto{\Rrightarrow}
\def\PDer{\cPDER}

\newlength{\seplen}
\newlength{\sepwid}
\def\firstblank{\,\rule{\seplen}{\sepwid}\,}
\def\secondblank{\firstblank\llap{\raisebox{2pt}{\firstblank}}}

\providecommand{\abbrv}[1]{#1.\@\xspace}
  \providecommand{\ie}{\abbrv{i.e}}

\def\To{\Rightarrow}

\def\twoar{\mathchoice
{\ooalign{%
\raise-.1em\hbox{$\rightharpoondown$}\cr%
\raise.1em\hbox{$\rightharpoonup$}}}%
{\ooalign{%
\raise-.1em\hbox{$\rightharpoondown$}\cr%
\raise.1em\hbox{$\rightharpoonup$}}}%
{\ooalign{%
\raise-.08em\hbox{\scriptsize$\rightharpoondown$}\cr%
\raise.08em\hbox{\scriptsize$\rightharpoonup$}}}%
{\ooalign{%
\raise-.075em\hbox{\tiny$\rightharpoondown$}\cr%
\raise.075em\hbox{\tiny$\rightharpoonup$}}}%
}
\def\eq{\textsc{eq}}
\def\coeq{\textsc{coeq}}
\def\pt{\text{pt}}
\def\yon{\textsc{y}}

\usepackage{lmodern}

\DeclareFontEncoding{LS1}{}{}
\DeclareFontSubstitution{LS1}{stix}{m}{n}
\DeclareSymbolFont{symbols4}{LS1}{stixbb}{m}{it}
\SetSymbolFont{symbols4}{bold}{LS1}{stixbb}{b}{it}
\DeclareMathSymbol{\rightmoon}{\mathbin}{symbols4}{"F9}
\DeclareMathSymbol{\leftmoon}{\mathbin}{symbols4}{"FA}

\usepackage{turnstile}
  \newcommand{\adjunct}[2]{\nsststile{#2}{#1}}

\makeatletter
\providecommand{\leftsquigarrow}{
  \mathrel{\mathpalette\reflect@squig\relax}
}
\newcommand{\reflect@squig}[2]{
  \reflectbox{$\m@th#1\rightsquigarrow$}
}
\makeatother

\newcommand{\var}[2]{ \begin{smallmatrix} #1 \\ \downarrow \\ #2 \end{smallmatrix}}
\newcommand{\smat}[1]{\left[ \begin{smallmatrix} #1 \end{smallmatrix}\right]}

\def\Ran{\text{Ran}}

\def\lan{\text{lan}} 
\def\Lan{\text{Lan}}

\newcommand{\mvs}{%
\fontfamily{mvs}\fontencoding{U}\fontseries{m}\fontshape{n}\selectfont}
\newcommand{\mvchr}[1]{{\mvs\char#1}}

\def\leaden{\text{\mvchr{198}}}
\def\iff{\Leftrightarrow}

\def\chor{\mathbb{Chr}}
\def\dpfs{\textsc{dpfs}\@\xspace}
\def\dfs{\textsc{dfs}\@\xspace}

\usepackage{tikz,tikz-cd}
\usetikzlibrary{arrows}
\usetikzlibrary{babel}

\tikzset{
	implies/.style={%
		 double
		,double equal sign distance
		,-implies
	},
	shorten <>/.style={
		 shorten >=#1
		,shorten <=#1}
	}

\newlength{\spcng}
\def\scaling{.2}
\newlength{\rsng}

\def\drawluangle{\begin{tikzpicture}[scale=\scaling]
\draw (0,0) -- (0,1) -- (1,1);%
\end{tikzpicture}}
\def\drawrdangle{\begin{tikzpicture}[scale=\scaling]
\draw (0,0) -- (1,0) -- (1,1);%
\end{tikzpicture}}
\def\rdangle   {\raisebox{\rsng}{\drawrdangle}}
\def\luangle   {\raisebox{\rsng}{\drawluangle}}

\def\boite{\square}

\newcommand{\prepull}{\protect{\mathchoice{\rdangle}{\rdangle}{\lrcorner}{\lrcorner}}}

\tikzset{
    ultra thin/.style= {line width=0.1pt},
    very thin/.style=  {line width=0.2pt},
    thin/.style=       {line width=0.35pt},% thin is the default
    semithick/.style=  {line width=0.6pt},
    thick/.style=      {line width=0.8pt},
    very thick/.style= {line width=1.2pt},
    ultra thick/.style={line width=1.6pt}
}
\def\cdisplayed{\mathbin{\ooalign{\hfil\raisebox{.225em}{$\shortmid$}\hfil\cr\raisebox{-.12em}{\scalebox{1}[1]{$=$}}\cr}}}
\def\cnormaled{\cdisplayed}
\def\cscripted{\ooalign{\hfil\raisebox{.14em}{$\scriptstyle\shortmid$}\hfil\cr\raisebox{-.12em}{\scalebox{1}[1]{$\scriptstyle =$}}\cr}}
\def\cscriptscripted{\ooalign{\hfil\raisebox{.05em}{$\scriptscriptstyle\shortmid$}\hfil\cr\raisebox{-.12em}{\scalebox{1}[1]{$\scriptscriptstyle =$}}\cr}}

\def\corth{\mathchoice{\cdisplayed}{\cnormaled}{\cscripted}{\cscriptscripted}}%

\def\lcorth#1{\leftidx{^{\corth}}{}{#1}}
\def\prescript#1#2#3{\leftidx{^{#1}}{_{#2}}{#3}}
\def\lorth#1{\prescript{\perp}{}{#1}}

\usepackage{xparse}
\usepackage{expl3}

\ExplSyntaxOn

\cs_new:Nn \if_intersects:nnTF {
  \tl_set:Nn \l_tmpa_tl {#1}
  \bool_set_false:N \bool_tmpa_l
  \tl_map_inline:nn {#2} {
    \tl_if_in:NnTF \l_tmpa_tl {##1} {
      \tl_map_break:n { \bool_set_true:N \bool_tmpa_l }
    } {}
  }
  \bool_if:NTF \bool_tmpa_l { #3 } { #4 }
}

\NewDocumentCommand{\IfIntersects}{mmmm}{
  \if_intersects:nnTF {#1} {#2} {#3} {#4}
}

\ExplSyntaxOff

\usepackage{tikz}

\tikzset{
  maybedashed/.code n args={2}{\IfIntersects{#1}{#2}{\pgfqkeys{/tikz}{dashed}}{\pgfqkeys{/tikz}{solid}}}
}

\NewDocumentCommand{\DrawCube}{m}{
  \begin{tikzpicture}[scale=.165]
    \coordinate (BLD) at (0,0);
    \coordinate (BLU) at (0,1);
    \coordinate (BDR) at (1,0);
    \coordinate (BRU) at (1,1);
    \coordinate (LDF) at (0-.5,0-.5);
    \coordinate (LUF) at (0-.5,1-.5);
    \coordinate (DRF) at (1-.5,0-.5);
    \coordinate (RUF) at (1-.5,1-.5);

    \IfIntersects{#1}{B}{ \fill[opacity=0.025] (BLD) -- (BLU) -- (BRU) -- (BDR); }{}
    \IfIntersects{#1}{L}{ \fill[opacity=0.025] (BLD) -- (BLU) -- (LUF) -- (LDF); }{}
    \IfIntersects{#1}{D}{ \fill[opacity=0.025] (BLD) -- (BDR) -- (DRF) -- (LDF); }{}
    \IfIntersects{#1}{R}{ \fill[opacity=0.025] (BDR) -- (BRU) -- (RUF) -- (DRF); }{}
    \IfIntersects{#1}{U}{ \fill[opacity=0.025] (BLU) -- (BRU) -- (RUF) -- (LUF); }{}
    \IfIntersects{#1}{F}{ \fill[opacity=0.025] (LDF) -- (LUF) -- (RUF) -- (DRF); }{}

    \IfIntersects{#1}{BL}{ \draw[thin] (BLU |- LUF) edge[maybedashed={#1}{F}] (BLD) edge[maybedashed={#1}{U}] (BLU); }{}
    \IfIntersects{#1}{BD}{ \draw[very thin] (BDR -| RUF) edge[maybedashed={#1}{F}] (BLD) edge[maybedashed={#1}{R}] (BDR); }{}
    \IfIntersects{#1}{LD}{ \draw[very thin] (BLD) edge[maybedashed={#1}{F}] (LDF); }{}
    \IfIntersects{#1}{BU}{ \draw[very thin] (BLU) -- (BRU); }{}
    \IfIntersects{#1}{LU}{ \draw[very thin] (BLU) -- (LUF); }{}
    \IfIntersects{#1}{BR}{ \draw[very thin] (BDR) -- (BRU); }{}
    \IfIntersects{#1}{DR}{ \draw[very thin] (BDR) -- (DRF); }{}
    \IfIntersects{#1}{RU}{ \draw[very thin] (BRU) -- (RUF); }{}
    \IfIntersects{#1}{LF}{ \draw[very thin] (LDF) -- (LUF); }{}
    \IfIntersects{#1}{DF}{ \draw[very thin] (LDF) -- (DRF); }{}
    \IfIntersects{#1}{UF}{ \draw[very thin] (LUF) -- (RUF); }{}
    \IfIntersects{#1}{RF}{ \draw[very thin] (DRF) -- (RUF); }{}
  \end{tikzpicture}
}

\usepackage{hyphenat}
\setlist[1]{itemsep=0pt}

\title{Localization theory for derivators}
\author{Fosco Loregian}
\date{\today}
\address{%
	\textsf{Fosco Loregian}: \newline
	Masaryk University, Faculty of Sciences              \newline
	Kotl\'{a}\v{r}sk\'{a} 2, 611 37 Brno, Czech Republic \newline
	\href{mailto:loregianf@math.muni.cz}
     {\sf loregianf@math.muni.cz}}
\allowdisplaybreaks

\begin{document}

\begin{abstract}
We outline the theory of reflections for prederivators, derivators and stable derivators. In order to parallel the classical theory valid for categories, we outline how reflections can be equivalently described as categories of fractions, reflective factorization systems, and categories of algebras for idempotent monads. This is a further development of the theory of monads and factorization systems for derivators. %We explore the definition of monadic adjunction between derivators.
\end{abstract}

\maketitle

%!TEX root = ../localization.tex

{\scriptsize \tableofcontents}

\section{Introduction}
\label{sec:intro}
The notion of \emph{co/reflection} (or \emph{co/reflective localization}) $i : \cB \rightleftarrows \cC$ combines together two of the most natural and pervasive notions in category theory: it is an adjunction between $\cC$ and one of its full subcategories $\cB$.

Since (echoing \cite{McL}) `adjoints are everywhere', it is easy to believe that reflections pop up quite often in Mathematics: indeed, many theorems inside and outside category theory admit nifty translations in terms of the existence of a certain co/reflective localization. Even more, certain theorems are all about describing certain classes of well-behaved categories in terms of reflective localizations of `prototypical' such categories: for example, a fundamental characterization of \emph{locally finitely presentable} categories \cite{Adamek1994} is that they all arise as localizations of presheaf categories $[\A\op, \Set]$ on small categories $\A$, whereas \emph{Giraud theorem} refines this statement asserting that \emph{Grothendieck toposes} all arise as \emph{left-exact} localizations of such presheaf categories (this means that the left adjoint preserves finite limits). Moreover, we currently have a great deal of ways to characterize reflective subcategories of a given $\cC$ in terms of other data or structures on it: more precisely,
\begin{enumerate}[label=\textsc{l}\oldstylenums{\arabic*}), ref=\textsc{l}\oldstylenums{\arabic*}]
	\item \label{locitem1} co/reflective subcategories of $\cC$ arise as special \emph{categories of fractions}, \ie categories where a given class of arrows in $\cC$ has been formally inverted into a category $\cC[\cS^{-1}]$, initial with this property. More in detail, a co/reflective subcategory arises inverting all the arrows that the left adjoint sends to isomorphisms. Somehow understandably, homotopy theorists feel comfortable with the nomenclature ``reflective localization'' because these are particular cases of homotopy categories;
	\item \label{locitem2} co/reflective subcategories of $\cC$ are in bijection with (categories of co/algebras for) idempotent co/monads on $\cC$ \cite[§4.2]{Bor2} (recall that a monad $T$ is \emph{idempotent} if its multiplication $\mu : T^2\To T$ is invertible);
	\item \label{locitem3} co/reflective subcategories of $\cC$ are in bijection with \emph{co/reflective prefactorization systems} \cite{CHK} on $\cC$ (recall that a prefactorization system is \emph{reflective} if its right class satisfies the `two out of three' property).
\end{enumerate}
Somehow, this series of equivalences draws the state of the art on the subject: our aim in the present paper is to provide a similar description of localization theory in a the framework of pre/derivators.

The current explosive development of higher category theory forces (not only, but especially) the community of category theorists, geometers and topologists to re-enact many classical statements into the language of $(\infty,1)$-category theory, aiming to provide a robust framework in which to find newer and newer applications for higher category theory.

Now, approaching this subject everybody notices quite quickly that there are many models we can choose to work in: among many, we cite simplicial categories \cite{bergner2007model}, model categories \cite{Qui}, quasicategories \cite{Joy,HTT}; of course, choosing the `right' framework to address a particular problem is a matter of taste, mathematical experience, and --to a certain extent-- fashion. What is even more clear though is that category theorists are quite interested in comparing the same construction in different models, and (even better) having ways to prove their mutual equivalence. Even though there are formal ways to address this issue, there is no general recipe to build these dictionaries between models. As a result, it can be easy, difficult or extremely difficult to define the same notion, according to the framework we chose to work within, let alone to compare it to others.

The theory of co/reflective localizations makes no exception in this respect: in the setting of simplicial categories it is quite easy to rephrase what is a reflection functor, with the only care that not all unenriched functors $F : |\underline{\cA}|\to |\underline{\cB}|$ between the underlying categories of two simplicial categories lift to simplicial functors $F : \underline{\cA} \to \underline{\cB}$; in the setting of model categories the `homotopy invariant' notion of co/reflective localization is that of a left/right \emph{Bousfield localization} \cite{Hirschhorn2003}, and this can be seen as solving the universal problem of enlarging the class of weak equivalences $\cW$ of a model category $\cM$ into $\cW\subset \cW'$, while maintaining the co/fibrations fixed. This results in the homotopy category $\ho(\cM,\cW)$ being a localization (in the 1-categorical sense) of $\ho(\cM,\cW')$.\fxnote{Or vice-versa? Check} Finally, in the realm of quasicategories the process of localization is captured by the following construction: if $W\subset \cC$ determines a class of edges of the quasicategory $\cC$, then the \emph{localization} $\cC[W^{-1}]$ of $\cC$ at $W$ is defined by the homotopy pushout
\[
\xymatrix{
W \ar@{}[dr]|{\luangle}\ar[r]\ar[d]& \cC \ar@{.>}[d]\\
\widetilde{W} \ar@{.>}[r]& \cC[W^{-1}]
}
\]
in the category of simplicial sets, where $\widetilde{W}$ is the fibrant replacement of the simplicial subset $W\subset \cC$ in the Kan-Quillen model structure (a model for which is $\text{Ex}^\infty(W)$). 

Now, the theory of derivators, initiated in order to grasp a more intrinsic description of the construction exhibiting the derived category of an abelian category $\A$, is able to get rid of simplicial machineries in $(\infty,1)$-category theory (or at least, it reduces them to the bare minimum with which a category theorist or an algebraist is comfortable). Even though there is no clear evidence that the 2-category of prederivators, as defined in \cite{Grothendieck}, but especially in \cite{groth2013derivators} can be made (weakly) equivalent to one of the other models for $(\infty,1)$-categories (but see \cite{carlson:prederivators} for a partial result in this direction), there is a certain effort to establish to which extent this comparison of models is possible, were it only because derivators are `friendlier' in that they only appeal 2-dimensional category theory. The present paper is part of this effort: here, we study the theory of co/reflections for prederivators, and we provide equivalent characterizations for these co/reflections in terms of objects of fractions, algebras for idempotent monads, and reflective factorization systems echoing characterizations \ref{locitem1}---\ref{locitem3} above.
\subsection*{Organization of the paper.}
In \autoref{sec:review} we briefly review the fundamental definitions of 0-, 1- and 2-cells in the 2-category of prederivator, and the subsequent refinement to the full sub-2-category of \emph{derivators}, as well as the notions of co/sieve, homotopy exact and $\D$-exact square, etc.; this is meant to be a quick reference for the reader, but it may appear terse if approached without a previous knowledge of these definitions. 
In \autoref{sec:refl-loc} we provide the basic definition of a reflection in $\PDer$; it is precisely what we expect it to be in a generic 2-category: an adjunction $L\adjunct{\eta}{\varepsilon}R$ with invertible counit $\varepsilon$ (co-reflections are defined dually, of course). Here, in \autoref{fractions} we establish the equivalence between (left exact) reflections of $\D$ and derivators of fractions with respect to sub-prederivators $\sS\subset \D$ (admitting a left calculus of fractions), and their description (\autoref{refl_are_ortho_classes}) in terms of \emph{coherently orthogonal} classes as defined in \autoref{cohorthogo} (but more extensively used in \cite{tderiv}). Animated by the desire to get rid of the stability assumption in our previous work, we had to refine our notion of factorization system: this yield the notion of \emph{choric} factorization \autoref{choric_fs}, and its equivalence with \emph{choric reflections} (\autoref{choric_ref}). %In \autoref{sec:preslim} we deal with another cornerstone of the theory, preservation, reflection and creation of co/limits. 
In \autoref{sec:monads} we investigate the connection between co/localizations and co/algebras for idempotent co/monads: the equivalence result is easily reached (\autoref{char_of_ref_using_mnds}), but leaves us with the feeling that outlining such a pervasive and useful theory as that of monads is an urgent matter: we take advantage of the theory initiated in \cite{lagkas} expanding it and polishing some of its corners, and we prove in \autoref{bij_with_ref} the equivalence between left exact localization and algebras for left exact idempotent monads.%; we attempt to state the definition of monadic morphism of derivators \autoref{monadic}: classically, a milestone for such a theory is the Barr-Beck monadicity theorem, that says an adjunction $F\dashv G$ is monadic if and only if $G$ is conservative and creates certain shapes of colimits; hopefully, subsequent work will be devoted to investigate the possibility to prove a similar result for derivators.

\subsubsection*{Addendum.} A more general version of some results presented here has been obtained independently by I. Coley \cite{coley} (see §3.3.3, §4.4.10, the account of regularity property, and 7.14, that works as a motivation for that paper). We remain available to the author and the readers of both papers, welcoming the notification of any other potential overlap to properly acknowledge it.
\subsection*{Notation and terminology.}\leavevmode
\subsubsection*{Foundations.}	
Among different foundational conventions that one may adopt, in this paper we assume that every set lies in a suitable Gro\-then\-dieck universe \cite{artin1972sga}. This choice can nevertheless be safely replaced by the more popular (albeit less powerful) foundation using sets and classes.

More in detail we implicitly fix a universe $\mho$, whose elements are termed \emph{sets}; \emph{small categories} have a \emph{set} of morphisms; \emph{locally small} categories are always considered to be small with respect to \emph{some} universe: treating with derivators it is a common choice to employ the so\hyp{}called \emph{two\hyp{}universe convention}, where we postulate the existence of a universe $\mho^+\succ \mho$ in which all the classes of objects of non\hyp{}$\mho$\hyp{}small, locally small categories live. 

\subsubsection*{Categories and functors.} Possibly large categories will be usually denoted as calligraphic letters like $\cA,\cB,\cC$ and suchlike; classes of morphisms in a category, often confused with the subcategory they generate, are denoted as calligraphic letters $\cE, \cM, \cX,\cY,\dots$ as well; when they are considered as objects of the category of categories, small categories are usually denoted as capital Latin letters like $I,J,K$ and suchlike: we denote in the same way an object of a possibly large category $\cC$; this slight abuse of notation causes no harm whatsoever. The 2-categories of diagrams, small categories, categories, prederivators and derivators, and more generally all 2-categories, are denoted in a sans\hyp{}serif typeface like $\Dia,\Cat,\CAT,\PDer,\Der$. The correspondence that inverts 1-cells of a 2-category is denoted $(\firstblank)\op$, whereas the correspondence that inverts 2-cells is denoted $(\firstblank)^\text{co}$. We denote $J_{/j}$ the \emph{slice} category of $J$ at the object $j$, and having objects the arrows with codomain $j$; dually, $J_{j/}$ denotes the category of morphisms having domain $j$.

Functors between \emph{small} categories are usually denoted as lowercase Latin letters like $u,v,w,\dots$ and suchlike (there must be of course numerous deviations to this rule); an hom\hyp{}object in a category $\cK$ or higher category $\mathsf{K}$ is often denoted $\cK(A,B)$ or $\mathsf{K}(A,B)$: a category of functors makes no exception, so we denote it $\Cat(A,B)$ and $\CAT(\cA,\cB)$; the symbols $\firstblank$, $\secondblank$ are used as placeholders for the ``generic argument'' of any kind of mapping; natural transformations between functors; or more generally 2-cells in a 2-category, are often written in Greek, or Latin lowercase alphabet, and collected in the set $\Nat(F,G)$.

Whenever there is an adjunction $F\dashv G$ between functors, that we denote $F:\cA \leftrightarrows \cB : G$, the arrow $Fa\to b$ in the codomain of $F$ and the corresponding arrow $a\to Gb$ in its domain are called \emph{mates} or \emph{adjuncts}; so, the notation ``the mate/adjunct of $f\colon Fa\to b$'' means ``the unique arrow $g\colon a\to Gb$ determined by $f\colon Fa\to b$''. When there is an adjunction between two functors $F,G$ we adopt $F\adjunct{\eta}{\epsilon}G$ as a compact notation to denote all at once that $F$ is left adjoint to $G$, with unit $\eta \colon 1 \to GF$ and counit $\epsilon\colon FG\to 1$. A customary choice of notation for the \emph{whiskering} between a 1-cell $F$ and a 2-cell $\alpha$ is $F * \alpha$ or $\alpha * F$. In order to avoid confusion with the many occurrences of an `upper-star' besides a morphism, we choose to denote the pre- and post-composition morphisms induced by $F : \D \to \E$ as $F_\dag : \PDer(\sX,\D) \to \PDer(\sX, \E)$ and $F^\dag : \PDer(\E,\sX) \to \PDer(\D, \sX)$ respectively.

\subsubsection*{Special categories.} Derivator theory forces to work with a huge variety of category shapes, and forces to choose clever notation to denote these categories $I,J,\dots$ as well as the functors $I\leftrightarrows J$. In our work, many choices are classical: for example, the simplex category $\bDelta$ is the \emph{topologist's delta}, having objects \emph{nonempty} finite ordinals $[n]=\{0 < 1 < \dots < n\}$ (this is opposed to the \emph{algebraist's delta} $\bDelta_+$ which has an additional initial object $[-1]$); we denote $\Delta^n$ the representable presheaf on $[n]\in\bDelta$, \ie the image of $[n]$ under the Yoneda embedding of $\bDelta$ in the category of simplicial sets $\widehat{\bDelta}$. More often though, the objects of $\bDelta$ are considered as categories via the obvious embedding: as a consequence, certain objects have many names (for example, the terminal object $[0]$ of $\bDelta$ and $\Cat$ is called $e$ in $\Dia$). 

The notation for other common categories deserves to be explained; the ``generic span'' $\{2\leftarrow 0\to 1\}$ will be denoted as $\luangle\,$, where the opposite category ``generic cospan'' $\{0\to 2\leftarrow 1\}$ will be denoted as $\prepull\,$. The \emph{nerves} of these two categories are the simplicial sets $\Lambda^2_0$ and $\Lambda^2_2$ (as it is customary to blur the distinction between a category and its nerve, we don't insist in keeping these notation separated). The completions of these two categories to ``generic commutative squares'' are obtained introducing a terminal (resp., initial) object into $\luangle$ (resp., $\prepull$), in such a way that these two categories have objects labeled $\begin{smallmatrix} 0&\to&1\\\downarrow&&\downarrow\\2&\to&\infty \end{smallmatrix}$ (resp., $\begin{smallmatrix} -\infty&\to&0\\\downarrow&&\downarrow\\1&\to&2 \end{smallmatrix}$); this choice permits agile notation as $X_{-\infty,0}, X_{-\infty,1}$ etc. to refer to the various sides of $X\in\D(\boite)$. Of course, these two categories are isomorphic, hence indistinguishable: it is only context that gives to their objects different labels. Another useful convention, employed from time to time to refer to the sides of $X\in\D(\boite)$ is the following: we write the sides $\begin{smallmatrix} X_{00} &\to& X_{10} \\ \downarrow && \downarrow \\ X_{01} &\to& X_{11}\end{smallmatrix}$ as $X_N$ (``$X$ north''), $X_S$ (``$X$ south''), etc., with evident meaning. The category $\twoar$ denotes the ``generic pair of parallel arrows'' $\{0\rightrightarrows 1\}$; its completion to a category with initial object is $\eq = \{-\infty\to0\rightrightarrows 1\}$, whereas its completion to a category with terminal object is $\coeq = \{0\rightrightarrows 1 \to +\infty\}$. The category $\bbtwo\times\bbtwo\times\bbtwo$, that looks like a cube, appears in certain arguments involving an adjunction ``lifted'' from $L\dashv R$ to $L^{\bbtwo}\dashv R^{\bbtwo}$ and restricted to a certain subcategory of the domain of $R^{\bbtwo}$; in these cases, we employ the \emph{build a roof} notation used in Rubik's cube theory that assigns to the faces of a cube letters {\tt BLDRUF} as in
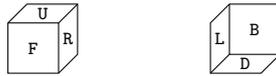
\begin{center}
\begin{tikzpicture}[scale=.666]
\draw 
  (0,0) coordinate (R) {} -- 
++(0,1) coordinate (A) {} --
++(45:\len) coordinate (P) {} --
++(1,0) coordinate (S) {} --
++(0,-1) coordinate (B) {} --
++(-135:\len) coordinate (T) {} -- cycle;
\draw (A) -- ++(1,0) -- (T); 
\draw (1,1) coordinate (Q) {} -- ++(45:\len);
\node at ($(P)!.5!(Q)$) {\tiny \tt U};
\node at ($(1,1)!.5!(B)$) {\tiny \tt R};
\node at ($(Q)!.5!(R)$) {\tiny \tt F};
\begin{scope}[xshift=4cm]
\draw 
  (0,0) coordinate (R) {} -- 
++(0,1) coordinate (A) {} --
++(45:\len) coordinate (P) {} --
++(1,0) coordinate (S) {} --
++(0,-1) coordinate (B) {} --
++(-135:\len) coordinate (T) {} -- cycle;
\draw (P) -- ++(0,-1) -- (B);
\draw (45:\len) -- ++(-135:\len);
\node at ($(P)!.5!(B)$) {\tiny \tt B};
\node at ($(A)!.5!(45:\len)$) {\tiny \tt L};
\node at ($(45:\len)!.5!(T)$) {\tiny \tt D};
\end{scope}
\end{tikzpicture}
\end{center}
according to the position of the face with respect to the observer. We choose to orient the cube in such a way there is a morphism of squares ${\tt F}\to {\tt B}$; when needed, we refer to the various edges of each face according to the aforementioned NSWE notation for $X\in\D(\boite)$. The subsets of the cube \DrawCube{BLDRUF} are then depicted \DrawCube{BD}, \DrawCube{BDF}, \DrawCube{BLR}, \DrawCube{LDF}, etc.
%!TEX root = ../localization.tex
\section{Review of derivator theory}
\label{sec:review}
\begin{defn}[Category of diagrams]
Let us start recalling that a {\em category of diagrams} is a full $2$-subcategory $\Dia$ of $\cCat$, such that
\begin{enumerate}[label=$\textsc{Dia}\oldstylenums{\arabic*})$]
\item all finite posets, considered as categories, belong to $\Dia$;
\item given $I\in \Dia$ and $i\in I$, the slice constructions $I_{i/}$ and $I_{/i}$ belong to $\Dia$;
\item if $I\in \Dia$, then $I\op\in \Dia$;
\item for every Grothendieck fibration $u : I\to J$, if all fibers $I_{j}$, for $j\in J$,  and the base $J$ belong to $\Dia$, then so does $I$. 
\end{enumerate}
\end{defn}
The minimal example for a category of diagrams is the locally posetal $2$-category $\textsf{fPos}$ of finite posets, while the maximal is $\Cat$, containing all small categories. There are other possible choices, like the $2$-category of finite categories.
\begin{defn}[the 2-category of prederivators]
If $\Dia$ is a category of diagrams, we call a \textbf{prederivator of type $\Dia$} a 2-functor $\D : \Dia\op \to \CAT$; the 2-category of prederivators has objects the prederivators, 1-cells the \emph{pseudonatural transformations}, and 2-cells the \emph{modifications} between pseudonatural transformations.\footnote{Taking \cite{groth2013derivators} as a standard reference for all the unexplained notation, we stick to the choice to call ``functors'' between derivators the 1-cells of $\PDer$, and ``natural transformations'' between morphisms of derivators its 2-cells.}
\end{defn}
The notion of derivator is a refinement of the notion of prederivator, motivated by the desire to provide a satisfactory axiomatization for triangulated categories that only appeals 2-categorical language. A \emph{derivator} is then a prederivator that satisfies the following additional conditions (we mimic the labeling convention of \cite{groth2013derivators}):
\begin{enumerate}[label=$\text{Der}\oldstylenums{\arabic*})$, ref=$(\text{Der}\oldstylenums{\arabic*})$]
	\item \label{derax:1} The functor $\D(I\sqcup J)\to \D(I)\times\D(J)$ obtained from the canonical inclusions $i_I : I\to I\sqcup J\leftarrow J : i_J$ is an equivalence.
	\item \label{derax:2} Each object $j : e\to J$ induces a family of functors $\D(J)\xto{j^*} \D(e)$; this family is jointly reflective, \ie a morphism $f\in \D(J)$ is invertible if and only if each $j^*f$ is invertible in $\D(e)$.
	\item \label{derax:3} Each functor $u^* : \D(J)\to \D(I)$ induced by $u : I\to J$ admits both a left adjoint $u_!$ and a right adjoint $u_*$. These functors are called, respectively, the homotopy left Kan extension and homotopy right Kan extension along $u$.
	\item \label{derax:4} Given a functor $u : J\to K$, there exist two squares in $\CAT$, induced by the colax pullbacks defining the slice and coslice categories, \ie by 2-cells in $\Dia$
\begin{align*}
		\vcenter{\xymatrix{
		J_{\!/k} \ar[r]^t\ar[d]_p & e\ar[d]^k\\
		J \ar[r]_u & K
		\ultwocell<\omit>{\varpi}
		}}
		\quad
		\rightsquigarrow
		\quad
		\vcenter{\xymatrix{
		\D(J_{\!/k})\ar[r]^{t_*} & \D(e)\dltwocell<\omit>{\varpi_*}\\
		\D(J)\ar[u]^{p^*} \ar[r]_{u_*} & \D(K)\ar[u]_{k^*}
		}}\\
		\vcenter{\xymatrix{
		J_{k/} \drtwocell<\omit>{\varpi'} \ar[r]^t\ar[d]_p & e\ar[d]^k\\
		J \ar[r]_u & K
		% \ultwocell<\omit>{\varpi}
		}}
		\quad
		\rightsquigarrow
		\quad
		\vcenter{\xymatrix{
		\D(J_{k/})\ar[r]^{t_!} & \D(e)\\
		\D(J)\urtwocell<\omit>{\varpi_!}\ar[u]^{p^*} \ar[r]_{u_!} & \D(K).\ar[u]_{k^*}
		}}
	\end{align*}
	These squares are filled by invertible 2-cells $\varpi'_! : t_!p^*\To k^* u_!$, and $\varpi_* : k^*u_* \To t_* p^*$ for every derivator $\D$.
\end{enumerate}
\begin{rmk}
Taken all together, the axioms of derivator are meant to ensure that we can build a category theory which is expressive enough for concrete applications, and in particular, they are meant to express the fact that we can compute left and right Kan extensions for every functor $u : I\to J$ \ref{derax:3}, and that these extensions are pointwise \ref{derax:4}. 
\end{rmk}
\begin{rmk}
It is helpful to keep in mind the two paradigmatic examples of constructions giving rise to derivators:
\begin{itemize}
	\item Let $\cM$ be a model category; then the association $J\mapsto \ho(\cM^J)$ defines a derivator.
	\item Let $\cC$ be a quasicategory; then the association $J\mapsto \ho(\cC^{NJ})$ defines a prederivator, which is a derivator if $\cC$ is complete and cocomplete.
\end{itemize}
\end{rmk}
\begin{rmk}\label{the_coherence}
Even though we call them \emph{functors}, the 1-cells $F: \E\to \D$ of $\PDer$ are pseudonatural transformations between 2-functors; every such transformation comes equipped with \emph{coherence data} as part of its definition, and in particular with invertible 2-cells $\gamma_{F,u} : u^* \circ F_J \To F_I \circ u^*$, one for each $u : I\to J$, suitably compatible with composition.
\end{rmk}
\begin{rmk}\label{tehgammas}
Every such transformation induces additional 2-cells (the `adjuncts' of $\gamma_{F,u}$, or the 2-cells obtained by `base change' from $\gamma_{F,u}$: see \cite[§1.2]{groth2013derivators})
\[\footnotesize
\begin{array}{cc}
\xymatrix@R=1cm@C=1cm{
\E(I)\ar@{{}{ }{}}@/^.5pc/[dr]|{\epsilon_u\Swarrow} \ar[r]^{u_*}\ar@{=}@/_1pc/[dr]&\E(J)\ar@{}[dr]|{\gamma_F\Swarrow} \ar[r]^{F_J}\ar[d]|{u^*}&\ar@{{}{ }{}}@/_.5pc/[dr]|{\eta_u\Swarrow} \ar[d]|{u^*}\ar@{=}@/^1pc/[dr]\D(J)&\\
&\E(I)\ar[r]_{F_I}&\D(I)\ar[r]_{u_*}&\D(J)
}
&
\xymatrix@R=1cm@C=1cm{
\E(I)\ar@{{}{ }{}}@/^.5pc/[dr]|{\eta_u\Nearrow}\ar[r]^{u_!}\ar@{=}@/_1pc/[dr]&\E(J)\ar@{}[dr]|{\gamma^{-1}\Nearrow}\ar[r]^{F_J}\ar[d]|{u^*}&\D(J)\ar@{{}{ }{}}@/_.5pc/[dr]|{\epsilon_u\Nearrow}\ar[d]|{u^*}\ar@{=}@/^1pc/[dr]&\\
&\E(I)\ar[r]_{F_I}&\D(I)\ar[r]_{u_!}&\D(J)
}
\end{array}
\]
We denote these 2-cells obtained pasting the units and counits of the adjunctions $u_!\dashv u^*\dashv u_*$ with $\gamma_{F,u}$, respectively $(\gamma_{F,u})_*$ and $(\gamma_{F,u})_!$. Somewhat sloppily, these morphisms measure how far is $F$ from commuting with the Kan extensions $u_!,u_*$ (see \cite[{\S}2.2]{groth2013derivators})
\end{rmk}
\begin{rmk}
The 2-cells $\varpi_*,\varpi_!$ of axiom \ref{derax:4} can be regarded as the 2-cells obtained by base change from the 2-cell $\varpi$.
\end{rmk}
We recall a few classical definitions from \cite{groth2013derivators} that will be useful later:
\begin{defn}\label{sieve-cosieve}
Let $u\colon J\to K$ be a fully faithful functor which is injective on objects.
\begin{enumerate}[label=$\roman*$)]
\item The functor $u$ is called a \textbf{cosieve} if whenever we have a morphism $u(j)\to k$ in $K$ then $k$ lies in the image of $u.$
\item The functor $u$ is called a \textbf{sieve} if whenever we have a morphism $k\to u(j)$ in $K$ then $k$ lies in the image of $u.$
\end{enumerate}
\end{defn}
\begin{defn}
Let $\D$ be a derivator and let us consider a natural transformation $\alpha$ as indicated in the following square in $\Cat:$
$$
\xymatrix{
J_1\ar[r]^v\ar[d]_{u_1}\xtwocell[1,1]{}\omit &J_2\ar[d]^{u_2}\\
K_1\ar[r]_w&K_2
}
$$
The square is $\D$-\textbf{exact} if the base change $\alpha_!\colon{u_1}_!\circ v^\ast\to w^\ast\circ {u_2}_!$ (or, by \cite[Lemma 1.20]{groth2013derivators}, equivalently $\alpha_\ast\colon u_2^\ast\circ w_\ast\to v_\ast\circ u_1^\ast$) is a natural isomorphism. The square is called \textbf{homotopy exact} if it is~$\D$-exact for all derivators $\D$.
\end{defn}
\begin{rmk}
Note that axiom \ref{derax:4} can be rephrased as `the squares
\[
\xymatrix{
J_{\!/k} \ar[r]^t\ar[d]_p & e\ar[d]^k\\
J \ar[r]_u & K
\ultwocell<\omit>{\varpi}
}
\qquad
\xymatrix{
J_{k/} \drtwocell<\omit>{\varpi'} \ar[r]^t\ar[d]_p & e\ar[d]^k\\
J \ar[r]_u & K
}
\]
are homotopy exact.'
\end{rmk}
\begin{lem}\label{equiv_recognition}
Let $F : \D \to \E$ be a morphism of prederivators; then $F$ is an equivalence (\ie it admits an adjoint inverse $G$ such that $GF\cong 1_{\D}$ and $FG\cong 1_{\E}$) if and only if each component $F_J : \D(J)\to \E(J)$ is an equivalence of categories in $\CAT$.
\end{lem}
\begin{proof}
It is clear that the condition is necessary. To show that it is also sufficient, let $G_J$ be the adjoint inverse of $F_J$, in such a way that each $F_J\adjunct{\eta_J}{\varepsilon_J}G_J$ is an equivalence of categories; in particular, each unit $\eta_J$ and each counit $\varepsilon_J$ are invertible natural transformations. Given $u : I\to J$ now it's easily seen that the coherence morphisms $\gamma_{F,u} : u^* F_J \Rightarrow F_I u^*$, when pasted with $\eta,\varepsilon$, give isomorphisms
\[
\begin{tikzcd}
\E(J) \arrow[r, "G_J"] \arrow[rd, equal, bend right] & \D(J) \arrow[Rightarrow,ld,shorten <=.25em,shorten >=2.5em, "\varepsilon"] \arrow[r, "u^*"] \arrow[d, "F_J"] & \D(I) \arrow[Rightarrow,ld,shorten <=.2em,shorten >=1em, "\gamma_{F,u}"] \arrow[d, "F_I"] \arrow[rd, equal, bend left] & \arrow[Rightarrow,ld,shorten <=2em,shorten >=.5em, "\eta"] \\
{} & \E(J) \arrow[r, "u^*"'] & \E(I) \arrow[r, "G_I"'] & \D(I)
\end{tikzcd}
\] testifying that the components $G_J$ glue to a pseudonatural transformation $\E\to \D$.
\end{proof}
\begin{defn}[\cite{groth2013derivators}, 2.2]
Let $F :\D \to \E$ be a morphism of derivators and let $u : J\to K$ be a functor. The morphism $F$ \emph{preserves homotopy left} respectively \emph{homotopy right Kan extensions along} $u$ if the natural transformation
\[
(\gamma^F_u)_! : u_!F\to Fu_!\qquad\mbox{respectively}\qquad{\gamma^F_u}_\ast : Fu_\ast\to u_\ast F
\]
is an isomorphism.
\end{defn}
In order to discuss the interaction of limits and reflections more systematically, we include the following definition. These will turn useful in the proof of general facts about co/limits preservation properties of reflections.
\begin{defn}\label{defn:refl}
Let $F\colon\D\to\E$ be a morphism of derivators and let $u\colon A\to B$ be in $\cCat$.
\begin{enumerate}
\item The morphism $F$ \textbf{reflects left} (resp. \textbf{right}) \textbf{Kan extensions along $u$} if for every diagram $X\in\D(B)$ such that $FX\in\E(B)$ lies in the essential image of $u_!\colon\E(A)\to\E(B)$ (resp., of $u_*\colon\E(A)\to\E(B)$) then already $X$ lies in the essential image of $u_!\colon\D(A)\to\D(B)$ (resp., of $u_*\colon\D(A)\to\D(B)$).
\item The morphism $F$ \textbf{reflects colimits of shape $A$} if it reflects left Kan extensions along $A\to A^\rhd$ (resp., $A \to A^\lhd$).
\end{enumerate} 
\end{defn}
\begin{rmk}
\begin{enumerate}
\item A morphism $F$ of derivators reflects left Kan extensions along $u$ if and only if the opposite morphism $F\op$ reflects right Kan extensions along $u\op$.
\item Let $F\colon\cC\to\cD$ be a functor between complete and cocomplete categories. Then $F$ reflects colimits of shape $A$ if and only if $y_F$ preserves colimits of shape $A$. In fact, it suffices to note that $X\colon A^\rhd\to\cC$ is a colimiting cocone if and only if it lies in the essential image of the left Kan extension along $i_A\colon A\to A^\rhd$.
\end{enumerate}
\end{rmk}

\begin{rmk}\label{rmk:pres-refl-colim}
By \cite[Prop.~3.9]{groth:revisit} the following are equivalent for a morphism $F\colon\D\to\E$ of derivators and $A\in\cCat$.
\begin{enumerate}
\item The morphism $F$ preserves and reflects colimits of shape $A$. 
\item A cocone $X\in\D(A^\rhd)$ is colimiting if and only if the cocone $FX\in\E(A^\rhd)$ is colimiting.
\end{enumerate}
\end{rmk}

\begin{defn}\label{defn:cons}
A morphism of prederivators is \textbf{conservative} if every component is conservative.
\end{defn}

Thus, a morphism $F\colon\D\to\E$ is conservative if a morphism $f\colon X\to Y$ in $\D(A)$ is an isomorphism as soon as $Ff\colon FX\to FY$ is an isomorphism in $\E(A)$.

\begin{egs}\leavevmode\label{egs:ff-cons}
Fully faithful morphism of prederivators are conservative. In particular, reflections, co-reflections, and equivalences are conservative.
\end{egs} 
\begin{defn}\label{creates}
Let $F : \D \to \E$ be a morphism of prederivators; we say that $F$ \textbf{creates homotopy left} (\textbf{resp., right}) \textbf{Kan extensions along} $u : A\to B$ if for every $X\in \D(B)$, 
\begin{itemize}
	\item the object $X$ lies in the essential image of $u_!^\D$ (resp, $u_*^\D$) if and only if $F_BX$ lies in the essential image of $u_!^\E$ (resp, $u_*^\D$);
	\item $F$ preserves homotopy left (resp., right) Kan extensions along $u$.
\end{itemize}
\end{defn}
\begin{defn}
We say that a prederivator \textbf{has initial} (resp. \textbf{terminal}) \textbf{objects} if $t^* : \D(e)\to \D(\varnothing)$ (where $t : \varnothing \to e$ is the unique arrow) has a left (resp., right) adjoint. Of course the notion is meaningful when $\D$ is a derivator, as (see \cite[1.12.i]{groth2013derivators}) $\D(\varnothing)$ is the terminal category.
\end{defn}
Here, we record a rather useful characterization of morphisms of derivators commuting with finite homotopy limits.
\begin{lem}\label{megalemma}
The following conditions are equivalent for an adjunction $L : \sD \rightleftarrows \sE : R$ of derivators (the assumption that these are derivators can't be removed):
\begin{enumerate}[label=\textsc{cl}\oldstylenums{\arabic*}), ref=\textsc{cl}\oldstylenums{\arabic*}]
	\item \label{cl:uno}The left adjoint $L$ commutes with finite limits;
	\item \label{cl:due}The left adjoint $L$ commutes with terminal objects and homotopy pullbacks;
	\item \label{cl:ter}The left adjoint $L$ commutes with products and homotopy equalizers;
\end{enumerate}
\end{lem}
\begin{proof}
The proof that $\ref{cl:due}\Rightarrow \ref{cl:uno}$ is the dual of \cite[7.1]{Ponto2014}, where a non-elementary argument proved that one can actually construct homotopy colimits over $J\in\Dia$ out of coproducts and homotopy pushouts if $J$ is homotopy finite.

We show that \ref{cl:ter} implies \ref{cl:due}: given $X\in\D(\prepull)$ there is a diagram $\bar X\in \D(\twoar)$ such that $\holim_\prepull X\cong \holim_{\twoar} \bar X$.

To this end, define the following functors between finite categories:
\begin{itemize}
	\item $i_\prepull : \prepull \to \boite$ is the inclusion of the left and bottom edge of $\prepull\,$;
	\item the functor $i_\lhd$ adjoins an initial object to the category $\twoar$;
	\item the coequalizer $q_\prepull$ of the pair $d_0,d_1 : \bbone \rightrightarrows \prepull$, \ie the functor
	{\tiny \[ 
	\left[\vcenter{\xymatrix@R=3.5mm@C=4mm{
	& 0\ar[d]\\
	1 \ar[r] & 2 
	}}\right]
	\quad \xto{q} \quad
	\left[\vcenter{\xymatrix@C=4mm{
	(0\equiv 1)\ar@<2pt>[r]\ar@<-2pt>[r] & 2
	}}\right]
	\]}
	\item the coequalizer $q_\boite$ of the pair $(\firstblank)_N, (\firstblank)_W : \bbtwo \rightrightarrows \boite$, \ie the functor
	{\tiny \[\left[\vcenter{\xymatrix@R=3.5mm@C=4mm{
-\infty \ar[r]\ar[d]& 0\ar[d]\\
1 \ar[r] & 2 
}}\right]
\quad \xto{q_\square} \quad
\left[\vcenter{\xymatrix@C=4mm{
-\infty \ar[r]& (0\equiv 1)\ar@<2pt>[r]\ar@<-2pt>[r] & 2
}}\right]\]} 
\end{itemize}
Notice that these functors arrange in a commutative square so that $q_{\boite} i_{\prepull}=i_{\rightrightarrows}q_{\prepull}$, hence $(q_{\boite})_* (i_{\prepull})_*\cong(q_{\boite} i_{\prepull})_*=(i_{\rightrightarrows}q_{\prepull})_*\cong (i_{\rightrightarrows})_*(q_{\prepull})_*$. Thus, we are reduced to verify that $(-\infty)^*(q_{\boite})_*\cong (-\infty)^*$. Let $Y\in \mathbb D(\boite)$, to compute $(-\infty)^*(q_{\boite})_*Y$ one can use the axiom \ref{derax:4} (Kan extensions are pointwise in a derivator), so that
\[
(-\infty)^*(q_{\boite})_*Y\cong \holim_{(-\infty/q_{\boite})} \pr^*Y
\]
where $\pr : (-\infty/q_{\square})\to \square$ is the canonical functor from the comma category. It is easy to notice that the category $(-\infty/q_{\square})$ admits an initial object $\emptyset$. Hence, 
\[
\holim_{(-\infty/q_{\square})} \pr^*Y\cong \emptyset^* \pr^*Y\cong (-\infty)^*Y,
\]
since taking the homotopy limit over a category with an initial object $\emptyset$ is the same as evaluating at $\emptyset$. This concludes the proof, as this argument gives a canonical morphism
\[
L_{\twoar} q_{\prepull,*}X \xto{\quad} q_{\prepull,*}L_\prepull X
\] that can be easily seen to be invertible using \ref{derax:2}.\footnote{\cite[4.4.3.2]{HTT} gives a proof that homotopy pullbacks can be constructed using equalizers in the setting of quasicategories. The present proof can be seen as the analogue statement adapted to derivators (whenever $\D = \D_{\cC}$ for a quasicategory $\cC$, it is really equivalent): starting with $X\in \mathbb D(\prepull)$ we build a diagram whose incoherent image in $\mathbb D(e)^{\rightrightarrows}$ is
\[
X_0\times X_1\rightrightarrows X_2.
\]
This is exactly $\dia_{\twoar}(q_{\prepull})_*X$ constructed above.}
\end{proof}
%!TEX root = ../localization.tex
\section{Reflective sub\hyp{}prederivators}\label{sec:refl-loc}
\subsection{Generalities about reflections}
In this section we study some basics concerning co/reflections of derivators. Related to this topic, see also \cite{heller:htpythies,cisinski-tabuada:non-connective,tabuada:universal-invariants}. A rather precise study of reflection theory for derivators has also been given in \cite{coley} (see in particular §3.3.3 --that we proved independently as \autoref{prop:refl-loc}--, §4.4.10, the account of regularity property).

We begin by recalling the following notion from category theory.
\begin{defn}
An adjunction $L :\cC\rightleftarrows\cD : R$ in $\cCAT$ is a \textbf{reflection} if the right adjoint is fully faithful and it is a \textbf{co-reflection} if the left adjoint is fully faithful.
\end{defn}
Duality allows us to focus on reflections only. We will rarely mention co-reflections, as this dualization process is straightforward.
\begin{rmk}
Let $L :\cC\rightleftarrows\cD : R$ be a reflection in $\cCAT$.
\begin{enumerate}
\item If $\Sigma_L$ denotes the class of morphisms in $\cC$ which are inverted by $L$, then $\cC\to\cD$ is a model for the ``category of fractions'' $\cC\to\cC[\Sigma_L^{-1}]$ of $\cC$ at the class $\cS$ (see \cite[Prop.~5.3.1]{Bor1}), and this localization is reflective.
\item If $\cC$ is complete or cocomplete, then so is $\cD$ (\cite[Prop.~3.5.3, Prop.~3.5.4]{Bor1}).
\end{enumerate}
\end{rmk}
In the notation of the above remark, one can give a description of the essential image of the fully faithful right adjoint in terms of $\cS = \Sigma_L$, and we quickly recall the relevant notions.
\begin{defn}\label{elle-locals}
Let $L :\cC\rightleftarrows\cD : R$ be a reflection in $\cCAT$ and let $\cS$ be the class of morphisms inverted by $L$.
\begin{enumerate}
\item An object $X\in\cC$ is $\cS$-\textbf{local} if for every $f : Y\to Z$ in $\cS$ the precomposition 
\[
\xymatrix{
\cC(Z,X)\ar[rr]^{\firstblank\circ f}&&\cC(Y,X)
}
\]
is a bijection.
\item A morphism $f : Y\to Z$ in $\cC$ is an \textbf{$\cS$-local equivalence} if
\[
\xymatrix{
\cC(Z,X)\ar[rr]^{\firstblank\circ f}&& \cC(Y,X)
}
\]
is a bijection for every $\cS$-local object $X\in\cC$.
\end{enumerate}
\end{defn}
We denote the class of $\cS$-local objects in a category $\cC$ as $\cS_L$; sending a class of morphisms to its class of local objects, and then to the class of $\cS$-local equivalences sets up a closure operator $(\firstblank)^\text{lw}$ on $\hom(\cC)$; among the classes of $(\firstblank)^\text{lw}$-closed there are all $\Sigma_L$'s. More precisely,
\begin{lem}\label{lem:basic-refl-loc}
Let $L :\cC\rightleftarrows\cD : R$ be a reflection in $\cCAT$ and let $\cS_L$ be the class of morphisms inverted by $L$ (we shortly refer to these morphisms as the \emph{$L$-local equivalences}).
\begin{enumerate}
\item The counit $\varepsilon : LR\toiso\id$ is an isomorphism as is the image $L\eta : L\toiso LRL$ of the unit $\eta$ under $L$ (which is to say that $\eta_X\in \cS_L$ for $X\in\cD$).
\item An object $X\in\cC$ lies in the essential image of $R$ if and only if $\eta_X$ is invertible if and only if $X$ is $\cS_L$-local.
\item A morphism is an $\cS_L$-local equivalence if and only if it lies in $\cS$.
\item An object $X\in\cC$ is $\cS_L$-local if and only if 
\[
\xymatrix{
\cC(RLX,X)\ar[rr]^{\firstblank\circ \eta_X}&&\cC(X,X)
}
\]
is a bijection.
\end{enumerate}
\end{lem}
\begin{proof}
A right adjoint is fully faithful if and only if the counit is invertible. The invertibility of $L\eta$ follows from the triangular identities. Statements (2) and (3) are left as exercises. Finally, since $\eta_X$ lies by (1) in $\cS$, the map $\firstblank\circ\eta_X$ is bijective for every $\cS$-local $X$. Conversely, suppose that $\firstblank\circ\eta_X$ is bijective and let $\alpha : RLX\to X$ such that $\alpha\circ\eta_X=\id_X$. In order to check that $\eta_X\circ\alpha=\id_{RLX}$ it suffices to consider the bijection
\[
\xymatrix{
\cC(RLX,RLX)\ar[rr]^{\firstblank\circ\eta_X}&&\cC(X,RLX)
}
\]
for the $\cS$-local object $RLX$ and to note that $\eta_X\circ\alpha$ and $\id_{RLX}$ are both sent to $\eta_X$. 
\end{proof}
We are now ready to define the main notion we will handle throughout the paper: 
\begin{defn}\label{defn:refl-loc}
Let $L :\D\rightleftarrows\sE : R$ be an adjunction in $\PDer$ or $\cDER$.
\begin{enumerate}
\item The adjunction $L\adjunct{\eta}{\varepsilon} R$ is a \textbf{reflection} if $R$ is fully faithful, \ie the counit $\varepsilon : LR\toiso\id$ is invertible.
\item The adjunction $L\adjunct{\eta}{\varepsilon} R$ is a \textbf{co-reflection} if $L$ is fully faithful, \ie the unit $\eta :\id\toiso RL$ is invertible.
\end{enumerate}
\end{defn}
The notions of $\cS$-local objects and $\cS$-local equivalences can be extended in a straightforward way to this setting, to establish an analogue of \autoref{lem:basic-refl-loc}. In particular, the essential image of $R$ coincides with the class of $L$-local objects; this suggests the existence of an identification $\sE\cong \sD\llbracket\sS_L^{-1}\rrbracket$ for the domain of the right adjoint in a reflection, with respect to a suitable notion of ``object of fractions'' in $\PDer$. We establish such an equivalence in \autoref{prefractions} below.

If $L$ is a morphism of prederivators, we let $\sS_L$ denote the sub\hyp{}prederivator $\sS_L(J)$ of all $X \in\D^{\bbtwo}(J)$ such that $\dia_{\bbtwo} LX\in \D(J)^{\bbtwo}$ is an isomorphism. It is called the prederivator of $\sS$-\textbf{locals} (but our crossed reference remains \autoref{elle-locals}, where the notion is defined for categories, as we will use quite often the prederivator analogue of \autoref{lem:basic-refl-loc}).
\begin{defn}
Let $\sS \subseteq \sD^\bbtwo$ be a sub\hyp{}prederivator. We say that a morphism of prederivators $F  : \sD \to \sE$ \textbf{inverts} $\sS$ if for each $J\in\Dia$ and  $X\in \sS(J)$, $\dia_{\bbtwo}F^\bbtwo_J(X)$ is an isomorphism in $\sE(J)$.\footnote{Recall from \cite{groth2013derivators} that there is a natural functor $J\to [\D(J), \D(e)]$ regarding every morphism $\varphi : i\to j$ as a natural transformation of the corresponding classifying functors; in the cartesian closed structure of $\cCAT$ this corresponds to a functor $\D(J) \to [J, \D(e)]$, and we call this the \emph{underlying diagram} functor.}
\end{defn}
\begin{defn}[prederivators of fractions]\label{def:preder-of-fracs}
Let $\sD$ be a prederivator, and $\sS\subseteq \sD^{\bbtwo}$ a sub\hyp{}prederivator. Define $\PDer_{\sS}(\sD,\sE)$ to be the subcategory of $\PDer(\sD,\sE)$ made by all $F  : \sD \to \sE$ that invert $\sS$.

We call the \textbf{prederivator of fractions} of $\sD$ with respect to $\sS$ the prederivator (when it exists, unique up to isomorphism) $\sD\llbracket \sS^{-1}\rrbracket$ with a canonical morphism $\gamma  : \sD \to \sD\llbracket \sS^{-1}\rrbracket$ having the property that $L^\dag : \PDer(\D\llbracket \sS^{-1}\rrbracket,\sE)\to \PDer(\D,\sE)$ induces an equivalence of the domain with the sub-category $\PDer_{\sS}(\D,\sE)$, so that
\[
\PDer_{\sS}(\D,\sE) \cong \PDer(\D\llbracket \sS^{-1}\rrbracket,\sE).
\]
\end{defn}
It is worth to notice that the usual care is needed about the presence of a derivator of fractions $\D\llbracket \sS^{-1}\rrbracket$ in our universe. The standard choice amounts to choose a big enough universe of sets $\mathsf{SET}$ for which $\D\llbracket \sS^{-1}\rrbracket  : \Dia \to \CAT$ is a 2-functor.
\begin{prop}\label{prefractions}
The following conditions are equivalent for an adjunction $L : \sD \rightleftarrows \sE : R$ of prederivators:
\begin{itemize}
	\item $L\dashv R$ is a reflection of $\D$ onto $\E$;
	\item The pair $(\E, L)$ exhibits the prederivator of fractions $\D\llbracket \sS_L^{-1} \rrbracket$.
\end{itemize}
\end{prop}
\begin{proof}
The proof is a slick adaptation of a 2-categorical argument: let $L\adjunct{\eta}{\varepsilon} R$ be the adjunction in subject. Then $\PDer(R,\sX) \adjunct{\eta^\dag}{\varepsilon^\dag} \PDer(L,\sX)$ is an adjunction 
\[
\PDer(R,\sX) : \PDer(\sD,\sX) \rightleftarrows \PDer(\sE,\sX) : \PDer(L,\sX).
\]
which is moreover still a reflection (because $\varepsilon^\dag$ remains an isomorphism); the essential image of $L^\dag$ can be then easily characterized as the subcategory $\PDer_{\sS_L}(\sD,\sX)$.
\end{proof}
The notion of reflection can be specialized in several different ways: 
\begin{defn}\label{kinds_of_refle}
An adjunction of prederivators $L : \sD \rightleftarrows \sE : R$ is said to be
\begin{itemize}
\item \textbf{essential} if $L$ has a left adjoint $Z$;
\item \textbf{Frobenius} if it is essential and $Z\cong R$; 
\item \textbf{(left) exact} if $L$ commutes with finite right Kan extensions. 
\end{itemize}
\begin{defn}[regular prederivator]
Let $\D$ be a prederivator satisfying axiom \ref{derax:3}. We say that $\D$ is \textbf{regular} when given a finite category $K$ and when $\omega$ is regarded as a category in the obvious way, we have that the diagram
\[
\xymatrix@C=1.8cm{
\D(K\times\omega) 
\drtwocell<\omit>{\;\;\rho_{\D}}
\ar[r]^-{K\times\pt_{\omega,!}}
\ar[d]_{\holim_K\times\omega}& \D(K) \ar[d]^{\holim_K}\\
\D(\omega) \ar[r]_{\pt_{\omega,!}}& \D(e)
}
\]
commutes up to an invertible 2-cell $\rho_{\D}$.
\end{defn}
This definition comes from \cite[IV.5]{heller:htpythies}; see also \cite{cisinski-tabuada:non-connective,cisinski-tabuada:non-commutative}. Note that \cite{groth:char} every stable derivator is regular.
\begin{prop}
Let $\sD$ be a prederivator, and $L\dashv R : \sD \to \sE$ a \emph{left exact} reflection of $\sD$. If $\sD$ is regular, then so is $\sE$.
\end{prop}
\begin{proof}
The proof is simple, and follows analyzing the pasting diagram
\[
\xymatrix@C=1.8cm{
\sE(K\times\omega)
\ar[rrr]^-{K\times\pt_{\omega,!}}\ar[ddd]
&&& 
\sE(K)
\ar[ddd]\\
&
\D(K\times\omega)
\urtwocell<\omit>{\;\;\gamma}
\drtwocell<\omit>{\;\;\rho_{\D}}
\ar[ul]_{L_{K\times\omega}} 
\ar[r]^{K\times\pt_{\omega,!}}
\ar[d]_{\holim_K\times\omega}
&
\D(K)
\drtwocell<\omit>{\;\;\gamma'}
\ar[ur]_{L_K}
\ar[d]^{\holim_K}& \\
&
\D(\omega)
\drtwocell<\omit>{\;\;\gamma''}
\ar[dl]_{L_\omega}
\ar[r]_{\pt_{\omega,!}}
&
\D(e)
\ar[dr]_{L_e}& \\
\sE(\omega)
\uurtwocell<\omit>{<-4>\;\;\gamma'''}
\ar[rrr]_{\pt_{\omega,!}}
&&&
\sE(e)
}
\]
(the coincidence of this pasting square with the canonical filling $\rho_{\sE}$ is an easy check).
\end{proof}
The most important closure property enjoyed by a reflection of a prederivator in $\PDer$ is, however, under the properties that define  a derivator: this appears as \cite[3.3]{coley}, where the general theory of localizations is laid down in order to give a more modern account of the main theorem in \cite{Heller1997}.
\begin{prop}\label{prop:refl-loc}
Let $L :\D\rightleftarrows\E : R$ be an adjunction of prederivators.
\begin{enumerate}
\item If $\D$ is a derivator and $R$ is fully faithful, then also $\E$ is a derivator.
\item If $\E$ is a derivator and $L$ is fully faithful, then also $\D$ is a derivator.
\end{enumerate}
\end{prop}
\begin{proof}
By duality it suffices to take care of the first statement, and we establish each of the axioms \ref{derax:1}-\ref{derax:4} individually.
\begin{enumerate}[label=Der\oldstylenums{\arabic*})]
\item[\ref{derax:1}] Let $J_1$, $J_2\in \Dia$, and consider the obvious inclusions $J_1\xto{u_1}J_1\sqcup J_2\xleftarrow{u_2}J_2$. Consider now the following commuting diagram:
\[
\xymatrix{
\D(J_1\sqcup J_2)\ar[rr]^{[u_1^*,u_2^*]^t}\ar@<-2,5pt>[d]_{L_{J_1\sqcup J_2}}&&\D(J_1)\times\D(J_2)\ar@<-2,5pt>[d]_{L_{J_1}\times L_{J_2}}\\
\E(J_1\sqcup J_2)\ar[rr]^{[u_1^*,u_2^*]^t}\ar@<-2,5pt>[u]_{R_{J_1\sqcup J_2}}&&\E(J_1)\times\E(J_2)\ar@<-2,5pt>[u]_{R_{J_1}\times R_{J_2}}
} 
\]
By \ref{derax:1} for $\D$, the arrow $\D(J_1\sqcup J_2)\to \D(J_1)\times\D(J_2)$ is an equivalence. By commutativity, and using \autoref{lem:basic-refl-loc}, it is easy to check that this implies that the arrow $\E(J_1\sqcup J_2)\to \E(J_1)\times\E(J_2)$ is also an isomorphism. Furthermore, by \ref{derax:1} for $\D$, $\D(\emptyset)$ is not the empty category and so, by the existence of the functor $R_{\emptyset} : \D(\emptyset)\to \E(\emptyset)$, we deduce that $\E(\emptyset)$ is not the empty category.
\item[\ref{derax:2}] Let $f$ be a morphism in $\E(J)$ for some $J\in\Dia$. Then, $f$ is an isomorphism if and only if $R_J(f)$ is an isomorphism which happens, by \ref{derax:2}, if and only if $R_J(f)_j$ is an isomorphism for any $j\in J$. This last condition is clearly equivalent to the statement that $R_\bbone(f_j)$ is an isomorphism for any $j\in J$, which means exactly that $f_j$ is an isomorphism for any $j\in J$.
\item[\ref{derax:3}] Consider a functor $u : A\to B$ in $\Dia$ and let $u_\ast :\D(A)\to\D(B)$ be the right Kan extension in $\D$. For $X\in\E(A)$ we show that $u_\ast R_AX$ lies in the essential image of $R_B$. By \autoref{lem:basic-refl-loc} it suffices to show that the top horizontal morphism in
\[
\xymatrix{
\D(B)(R_BL_Bu_\ast R_AX,u_\ast R_AX)\ar[rr]^-{\eta u_*R_A}\ar[d]^\wr&&\D(B)(u_\ast R_AX,u_\ast R_AX)\ar[d]^\wr\\
\D(A)(u^\ast R_BL_Bu_\ast R_AX, R_AX)\ar[rr]^-{}\ar[d]^\wr&&\D(A)(u^\ast u_\ast R_AX,R_AX)\ar[d]^\wr\\
\E(A)(L_Au^\ast R_BL_Bu_\ast R_AX, X)\ar[rr]^-{}\ar[d]^\wr&&\E(A)(L_Bu^\ast u_\ast R_AX,X)\ar[d]^\wr\\
\E(A)(u^\ast L_BR_BL_Bu_\ast R_AX, X)\ar[rr]^\cong_-{u^*L_B\eta u_*R_A}&&\E(A)(u^\ast L_Bu_\ast R_AX,X)
}
\]
is an isomorphism. By naturality this diagram commutes, and the bottom horizontal is invertible by \autoref{lem:basic-refl-loc}. For every $X\in\E(B)$ and $Y\in\E(A)$ there is a chain of natural isomorphisms
\begin{align*}
\E(B)(X,L_Bu_\ast R_AY)\cong&\, \D(B)(R_BX,R_BL_Bu_\ast R_AY)\\
\cong&\, \D(B)(R_BX,u_\ast R_AY)\\
\cong&\, \D(A)(u^\ast R_BX,R_AY)\\
\cong&\, \D(A)(R_Au^\ast X,R_AY)\\
\cong&\, \E(A)(u^\ast X,Y),
\end{align*}
showing that $L_Bu_\ast R_A$ is a model for the right Kan extension along $u$ in $\E$. 

For the existence of left Kan extensions in $\E$, let $u_! :\D(A)\to\D(B)$ be a left Kan extension functor in $\D$, $X\in\E(A)$ and $Y\in\E(B)$. The natural bijections
\begin{align*}
\E(B)(L_Bu_!R_AX,Y)\cong&\, \D(B)(u_!R_AX,R_BY)\\
\cong&\, \D(A)(R_AX,u^\ast R_BY)\\
\cong&\, \D(A)(R_AX,R_Au^\ast Y)\\
\cong&\, \E(A)(X,u^\ast Y)
\end{align*}
show that $L_Bu_!R_A$ is a model for the left Kan extension along $u$ in $\E$.
\item[\ref{derax:4}] Let $u : A\to B$ be a functor in $\Dia$ and let $b\in B$. Using the explicit form of homotopy Kan extensions in $\E$ described above, it is not difficult to use the fact that homotopy Kan extensions are computed point-wise in $\D$ to show that the same happens in $\E$.\qedhere
\end{enumerate}
\end{proof}
\begin{rmk}\label{rmk:refl-loc}
Let $L :\D\rightleftarrows\E : R$ be a reflection of derivators and let $u : A\to B$. Left and right Kan extension functors in $\E$ along $u$ are respectively given by
\[
L_Bu_!R_A :\E(A)\to\E(B)\qquad\text{and}\qquad L_Bu_\ast R_A :\E(A)\to\E(B).
\]
\end{rmk}
\end{defn}
We record here the definition of a \emph{choric} reflection: sharpening \autoref{defn:refl-loc} in this way ensures that there is some control on the coherence with which the unit map $\eta : 1\To RL$ is given.
\begin{defn}[Choric reflection]\label{choric_ref}
Let $L\adjunct{\eta}{\varepsilon}R : \D \leftrightarrows \E$ be a reflection; we define 
$\chor_L$ to be the sub-prederivator %of $\D^{\bbtwo}$ of those $ X\in\D^{\bbtwo}$ such that $X_0$ lies in $\E(e)$:
\[
\chor_L  = \{ X\in\D^{\bbtwo} \mid X_1\in \E(e)\} \subseteq \D^{\bbtwo}
\]
Now, we define $R^\sharp : \chor_L\to \D$ to be the restriction of evaluation at $0$, \ie $R^\sharp \coloneqq 0^*|_{\chor_L}$.

We say that the reflection $(L,R)$ is \emph{choric} if $R^\sharp$ admits a left adjoint $L^\sharp$.
\end{defn}
\begin{rmk}
It is worth to spell out explicitly how some incoherent diagrams associated to the adjunction $L^\sharp\adjunct{\eta^\sharp}{\varepsilon^\sharp}R^\sharp$ look like. Every $X\in\D(\boite)$ induces $L^\sharp X \in \chor_L(\boite)\subseteq \D^{\bbtwo}(\boite) = \D(\DrawCube{BLDRUF})$, and
\[
\dia (X) = 
\vcenter{
{\tiny \xymatrix@R=.25cm@C=.25cm{
& RLX_{00} \ar[rr]\ar[dd]&& RLX_{01}\ar[dd]\\
X_{00} \ar[ur]\ar[rr]\ar[dd]&& X_{01}\ar[dd]\ar[ur]\\
& RLX_{10} \ar[rr]&& RLX_{11}\\
X_{10} \ar[ur]\ar[rr]&& X_{11}\ar[ur]
}}}
\]
in such a way the {\tt L} and {\tt R} faces of the cube correspond to naturality squares of the unit of $L\dashv R$, and the {\tt B} face correspond to the rule ``apply $L$ to the square $X$''.%\footnote{We employ the {\tt BLDRUF} notation for Rubik's cubes: the faces of a cube are denoted as {\tt U}p, {\tt D}own, {\tt L}eft, {\tt R}ight, {\tt F}ront, {\tt B}ack}
\end{rmk}
\begin{rmk}
As noticed, if $L\dashv R$ is a choric reflection, so that there is an induced adjunction $L^\sharp\dashv R^\sharp$, the latter adjunction gives a coherent choice of liftings of unit components for the former adjunction. This motivates the notation, as $(\firstblank)^\sharp$ ``lifts the adjunction (by a semitone)''. 

It may appear as if this notion is given with the only purpose to obtain an ad-hoc control of coherence for the unit $\eta$; it is nevertheless possible to prove that many reflections are choric (see \autoref{refl_are_fs}, where we prove that choric reflections arise as those associated to choric factorization systems (\autoref{choric_fs}), and every factorization system which is defined by an algebra structure for the squaring monad (\cite[Thm III]{tderiv}) on $\PDer$ is in fact choric).
\end{rmk}
We end this introductory subsection gathering a few examples of reflections of prederivators and derivators:
\begin{eg}
Let $i : \Dia \to \CAT$ be the inclusion functor. For every category $\cC$ there is the represented prederivator $y_\cC$ obtained via the map $\cC\mapsto \D_\cC  \coloneqq \CAT(i(\firstblank),\cC)$ (the ``restricted Yoneda embedding''), so that an adjunction $L :\cC\rightleftarrows\cD : R$ in $\cCAT$ is a reflection if and only if the induced adjunction $y_F : y_\cC\rightleftarrows y_\cD : y_G$ is.
\end{eg}
This has interesting specific sub-examples:
\begin{eg}
As a special case of particular relevance later there are the following adjunctions
\[
\cSp_{\geq 0}\rightleftarrows\cSp\qquad\text{and}\qquad\cSp\rightleftarrows\cSp_{\leq 0}
\]
exhibiting the derivator $\cSp_{\geq 0}$ of connective spectra as a co-reflection of $\cSp$ and the derivator $\cSp_{\leq 0}$ of coconnective spectra as a reflection of $\cSp$. More generally, let $k\in\lZ$. The full sub\hyp{}derivator $\cSp_{\leq k}\subseteq\cSp$ of $k$-coconnective spectra is closed under arbitrary  limits and the full sub\hyp{}derivator $\cSp_{\geq k}\subseteq\cSp$ of $k$-connective spectra is closed under  arbitrary colimits. This translates into a reflection between the associated derivators. 
\end{eg}
\begin{eg}
If $p$ is a prime number, the class of $p$\emph{-acyclic} chain complexes of abelian groups $\{X_\bullet\mid X_\bullet\otimes \mathbf{Z}/p\mathbf{Z}\simeq 0_\bullet\}$ is a co-reflection of $\text{Ch}(\mathbf Z)$. This translates into a reflection between the associated derivators.
\end{eg}
\begin{eg}
Let $F :\cM\rightleftarrows\cN : G$ be a Quillen adjunction between Quillen model categories and let
\[
\lL F  :\ho_\cM\rightleftarrows\ho_\cN : \lR G
\]
be the induced derived adjunction. If $(F,G)$ is a left Bousfield localization, then the derived adjunction is a reflection. Similarly, right Bousfield localizations induce co-reflections.
\end{eg}
\begin{eg}
The Kan extension adjunctions associated to fully faithful functors $u : A\to B$ yield co-reflections and reflections,
\[
u_!  :\D^A\rightleftarrows\D^B : u^\ast 
\qquad
\text{and}
\qquad 
u^\ast :\D^B\rightleftarrows\D^A : u_\ast.
\]
\end{eg}
\begin{eg}[Recollements of stable derivators]
Let $j : I\to J$ be a sieve (\autoref{sieve-cosieve}), and let $i : K \hookrightarrow J$ be the subcategory spanned by the complement of the image of $j$. Let $\D$ be a stable derivator. Then we can obtain recollements \cite{BBDPervers} of triangulated categories
\[
\xymatrix{
  \D(I) \ar[r]|{j_*} & \D(J)\ar[r]|{i^*}\ar@<8pt>[l]^{j^!}\ar@<-8pt>[l]_{j^*} & \D(K) \ar@<8pt>[l]^{i_*} \ar@<-8pt>[l]_{i_!}
& \D(K) \ar[r]|{i_!} & \D(J)\ar[r]|{i^*}\ar@<8pt>[l]^{i^*}\ar@<-8pt>[l]_{i^?} & \D(I)\ar@<8pt>[l]^{j_*} \ar@<-8pt>[l]_{j_!}
}\tag{$\star$}
\]
where the functors $j^!, i^?$ are respectively called co-exceptional and exceptional inverse images of $j$ and $i$.

This motivates the following definition: a \textbf{recollement of derivators} is an arrangement of morphisms of derivators
\[\xymatrix{
\D' \ar[r]|i & \ar@<-8pt>[l]_{i_L}\ar@<8pt>[l]^{i_R}\D\ar[r]|q  & \D''\ar@<-8pt>[l]_{q_L}\ar@<8pt>[l]^{q_R}
}\]
satisfying the following conditions
\begin{enumerate}[label=\textsc{r}\oldstylenums{\arabic*}), ref=\textsc{r}\oldstylenums{\arabic*}]
\item There are adjunctions $i_L\dashv i\dashv i_R$ and $q_L\dashv q\dashv q_R$;
\item The functors $i, q_L, q_R$ are all fully faithful;
\item The image of $i$ equals the \emph{essential kernel} of $q$, \ie the full subcategory of $\D$ such that $qX\cong 0$ in $\D^1$;
\item 
The natural homotopy commutative diagrams 
\[
\xymatrix{
  q_L q & \id_\D  & ii_R  & \id_\D \\
  0 & ii_L  & 0 & q_R q
  \ar "1,1";"1,2" ^{\epsilon_{(q_L\dashv q)}}
  \ar "1,1";"2,1" 
  \ar "1,2";"2,2" ^{\eta_{(i_L\dashv i)}}
  \ar "1,3";"1,4" ^{\epsilon_{(i\dashv i_R)}}
  \ar "1,3";"2,3" 
  \ar "1,4";"2,4" ^{\eta_{(q\dashv q_R)}}
  \ar "2,1";"2,2" 
  \ar "2,3";"2,4" 
}
\]
induced by the previous axioms are both cartesian and cocartesian.
\end{enumerate}
It is immediate to see that every sieve induces two recollements of derivators between the shifted derivators
\[
\xymatrix{
  \D^I \ar[r]|{j_*} & \D^J\ar[r]|{i^*}\ar@<8pt>[l]^{j^!}\ar@<-8pt>[l]_{j^*} & \D^K \ar@<8pt>[l]^{i_*} \ar@<-8pt>[l]_{i_!}
& \D^K \ar[r]|{i_!} & \D^J\ar[r]|{i^*}\ar@<8pt>[l]^{i^*}\ar@<-8pt>[l]_{i^?} & \D^I\ar@<8pt>[l]^{j_*} \ar@<-8pt>[l]_{j_!}
}
\]
such that the recollements in ($\star$) correspond to the evaluation of this diagram on the base. Recollements situation are quite natural ways to build reflections (and in fact more: a pair of reflective and bireflective sub\hyp{}prederivators) of a given stable derivator $\D$.
\end{eg}
More conceptual examples of reflections are the following:
\begin{eg}
Regarded as a 2-category, $\PDer$ supports a calculus of Kan extensions;\footnote{Every now and then we will implicitly rely on the following abstract result: $\PDer$ is at the same time the category of algebras for a 2-monad, and the 2-category of coalgebras for a 2-comonad, on the 2-category $[\Dia_0, \cCAT]$ ($\Dia_0$ is the class of objects of $\Dia$). More explicitly, there is a triple of adjoints
\[
\xymatrix{
	[\Dia\op,\Cat] \ar[r]|H & \ar@<7pt>[l]^{\Ran_H}\ar@<-7pt>[l]_{\Lan_H} [(\Dia\op)_0,\Cat]
}
\]
induced by the inclusion $H : \Dia_0 \hookrightarrow \Dia$. The monad of the adjunction $\Lan_H \dashv (\firstblank\circ H)$ and the comonad of the adjunction $(\firstblank\circ H)\dashv \Ran_H$ do the job.} in particular, we can define a 1-cell $G  : \D \to \sE$ in $\PDer$ to be \emph{dense} if the left extension $\lan_GG  : \sE \to \sE$ exists, is pointwise and isomorphic to the identity $\id_{\sE}$. 

In presence of a Yoneda structure on a 2-category $\mathsf{K}$, we can characterize a dense $G : \D \to \E$ as a 1-cell such that $\sE(G,1) = \lan_G\yon_{\D}$ is fully faithful; 
\[
\xymatrix{
\D \ar[r]^G \ar[d]_{\yon_{\D}} & \sE \ar@/^1pc/[dl]^{\lan_G\yon_{\D}}\\
\widehat{\D}
}
\]
when the 1-cell $\sE(G,1)$ admits a left adjoint, this determines a reflection of $\widehat{\D}$ (the \emph{Yoneda object} \cite{street1978yoneda} associated to $\D$). It is tempting to extend this characterization to the 2-category of prederivators (or a suitable sub-2-category thereof): a thorough discussion of all these issues will be the subject of a separate work.
\end{eg}

\subsection{Calculus of fractions}
Left exact reflections can be characterized via a derivator-theoretic analogue of \cite[5.6.1]{Bor1}: we first introduce the necessary terminology.
\begin{defn}%[Wide sub\hyp{}prederivator]
\label{wide}
Let $\sS\subseteq \sD^{\bbtwo}$ be a sub\hyp{}prederivator; we say that $\sS$ is \textbf{wide} if it is closed under composition\footnote{This means that whenever $X\in \D^\bbthree(J)$ is such that $X_{(0,1)}, X_{(1,2)} \in \sS(J)$ then also $X_{(0,2)}\in\sS(J)$; as a side note, we remark that in a similar fashion we can define \emph{cancellation properties} for a sub\hyp{}prederivator $\sS$: see \cite[3.17]{factonder} for more on this.} and contains all isomorphisms.
\end{defn}
\begin{defn}
\label{ore}
A sub\hyp{}prederivator $\sS\subseteq \sD^{\bbtwo}$ is said to \textbf{satisfy the right Ore condition} if for every $J\in\Dia$ the following two conditions are satisfied: 
\begin{enumerate}
	\item[o1)] for every $X\in \sD^J(\prepull)$ such that $X_{(1,2)}\in \sS(J)$, there exists an $X'\in \sD^J(\boite)$ such that $X'_{\prepull}\cong X$ and $X_{(-\infty,0)}\in \sS(J)$.
	\item[o2)] for every $X\in \sD^J(\coeq)$ such that $X_{(1,\infty)}$ is in $\sS(J)$, there exists an $X'\in \sD^J(\coeq)$ such that $X_{(-\infty,0)}$ is in $\sS(J)$.
\end{enumerate}
\end{defn}
\begin{defn}
A sub\hyp{}prederivator $\sS\hookrightarrow \sD^\bbtwo$ is said to \textbf{admit a right calculus of fractions} if it is wide and satisfies the right Ore condition.
\end{defn}
It turns out that choric reflection allow to find an analogue of \cite[5.6.1]{Bor1} for left exact reflections of derivators:
\begin{prop}\label{fractions}
The following conditions are equivalent for a choric reflection $L : \sD \rightleftarrows \sE : R$ of derivators:
\begin{enumerate}[label=\textsc{cf}\oldstylenums{\arabic*}), ref=\textsc{cf}\oldstylenums{\arabic*}]
\item \label{cf:uno}The left adjoint $L$ commutes with finite limits, so the reflection $L\dashv R$ is left exact (\autoref{kinds_of_refle});
\item \label{cf:due} The sub\hyp{}prederivator $\sS_L$ is \emph{homotopy pullback stable}: if $X\in\D(\boite)$ is cartesian and such that $X_E\in\sS_L(e)$, then $X_W \in \sS_L(e)$.
\item \label{cf:tre} The sub\hyp{}prederivator $\sS_L$ is closed under finite right Kan extensions in $\D$.
\end{enumerate}
\end{prop}
\begin{proof}
The proof proceeds in various steps that we state as different items: first, notice that \autoref{megalemma} above shows that \ref{cf:uno} is equivalent to the request that $L$ commutes with terminal objects and homotopy pullbacks, and to the request that $L$ commutes with products and homotopy equalizers. We will freely use this equivalence throughout the present proof.
\begin{enumerate}
\item To show that \ref{cf:uno} implies \ref{cf:tre} assume that $L\dashv R$ is a left exact reflection. Recall that \ref{cf:tre} says that whenever we have a diagram $Y\in \D(\DrawCube{BLDRUF})$ whose underlying diagram is
\[
\xymatrix@R=.25cm@C=.25cm{
& Y_{001} \ar[rr]\ar[dd]&& Y_{011}\ar[dd]^{s_3}\\
Y_{000} \ar[ur]\ar[rr]\ar[dd]_{s_1}&& Y_{010}\ar[dd]|(.35){s_2}\ar[ur]\\
& Y_{101} \ar[rr]&& Y_{111}\\
Y_{100} \ar[ur]\ar[rr]&& Y_{110}\ar[ur]
}
\]
where each vertical arrow $s_1,s_2,s_3$ is in $\sS_L$, and the {\tt U}, {\tt D} faces are cartesian, then also $s_4\in\sS_L$. Assuming \ref{cf:uno} and applying $L$ to it, we get that the {\tt U}, {\tt D} faces remain cartesian, and $Ls_4$ is the unique morphism connecting two pullbacks of isomorphic diagrams; this entails that $Ls_4$ is invertible.

\item We show that \ref{cf:due} implies \ref{cf:uno}. The functor $L$ commutes with finite right Kan extensions if and only if it commutes with finite limits \cite[2.4]{groth2013derivators}, and by our \autoref{megalemma}, this latter condition is true if and only if $L$ commutes with homotopy pullbacks (as $L$ already preserves terminal objects).

We show that given a cartesian square $i_{\prepull,*}Y = \begin{smallmatrix} Y_{-\infty} &\to& Y_0 \\ \downarrow && \downarrow \\ Y_1 &\to& Y_2\end{smallmatrix}$ where $Y_E$ lies in $\sS_L^{\corth}$, then $L$ preserves it, and then we show that this is in fact sufficient. Given such a $Y$, note that also $Y_E$ lies in $\sS_L^{\corth}$, as this class is closed under pullback. Now, we can embed $Y$ as the {\tt F} face of a cube $L^\sharp(Y)$ (as $L\dashv R$ is a choric reflection)
\[
{\scriptsize \xymatrix@R=.66em@C=.66em{
& RL Y_{-\infty} \ar[rr]\ar[dd]&& RL Y_0\ar[dd]\\
 Y_{-\infty} \ar[ur]\ar[rr]\ar[dd]&&  Y_0\ar[dd]\ar[ur]\\
& RL Y_1 \ar[rr]&& RL Y_2\\
 Y_1 \ar[ur]\ar[rr]&&  Y_2\ar[ur]
}}
\] Now $L$ commutes with pullbacks if and only if the {\tt B} face of this cube is cartesian; unwinding the definitions, and employing \ref{derax:2}, this is equivalent to say that in the diagram $L^\sharp Y \to i_{\boite,*}L^\sharp i_\boite^* Y$, whose underlying diagram is
\[
{\scriptsize 
\xymatrix@R=.25cm@C=.25cm{
& P \ar[rr]\ar[dd]&& RL Y_0\ar[dd]\\
RL Y_{-\infty} \ar[ur]^\zeta\ar[rr]\ar[dd]&& RL Y_0\ar[dd]\ar@{=}[ur]\\
& RL Y_1 \ar[rr]&& RL Y_2\\
RL Y_1 \ar@{=}[ur]\ar[rr]&& RL Y_2\ar@{=}[ur]
}}
\]
the arrow $\zeta  : LY_{-\infty} \to P$ is an isomorphism; this is equivalent to ask that $\zeta$ lies in $\sS_L\cap \sS^{\corth}_L$ (since $(\sS_L, \sS_L^{\corth})$ is a prefactorization). Now since $P\in\E(e)$, the arrow $\zeta$ lies in $\E(e)^{\bbtwo} = (\lorth{\E(e)^{\bbtwo})^\perp} = \sM$, so we only have to show that $\zeta$ gets inverted by $L$. To see this, consider the juxtaposition of squares
\[
\xymatrix{
Y_{-\infty} \ar[r]^{\eta_\infty}\ar[d] & RLY_{-\infty} \ar[r]^\zeta \ar[d] & P \ar[d] \\
Y_1 \ar[r]_{\eta_1} & RLY_1 \ar@{=}[r]& RLY_1;
}
\] the whole diagram is homotopy cartesian, so that (since we assume that $\sS_L$ is closed under pullback) $\zeta\circ \eta_\infty$ lies in $\sS_L$; but then, so do $\eta_1$ and $P \to RLY_1$, so (since $\sS_L$ is a 3-for-2 class), $\zeta\in\sS_L$.

It remains to show that now $L$ preserves pullbacks: given any $i_{\prepull,*}Y = \begin{smallmatrix} Y_{-\infty} &\to& Y_0 \\ \downarrow && \downarrow \\ Y_1 &\to& Y_2\end{smallmatrix}$ we can always $\fF$-factor its columns (because $(\sS_L,\sS_L^{\corth})$ is a \dfs and the functor $\D_{\fF}(\bbtwo) \to \D^{\bbtwo}(\bbtwo)$ is essentially surjective), to obtain
\[
{\scriptsize \xymatrix@R=.5cm@C=.5cm{
Y_{-\infty} \ar[r]\ar[d]_{s_\infty} &Y_0 \ar[d]^{s_0}\\
F \ar[r]\ar[d]_{s_\infty^\perp} & F'\ar[d]_{s_0^\perp}\\
Y_1 \ar[r] & Y_2 
}}
\]
now the upper square is cartesian (and the left column exhibits the factorization of $Y_W$) since $\sS_L$ is closed under pullback, so that $L$ preserves it; the lower square is cartesian, and the above argument shows that $L$ preserves it as well.

\item Now we prove that \ref{cf:tre} implies \ref{cf:due}: let $i_{\prepull,*}Y = \begin{smallmatrix} Y_{-\infty} &\to& Y_0 \\ \downarrow && \downarrow \\ Y_1 &\to& Y_2\end{smallmatrix}$ be a cartesian square such that $Y_E\in\sS_L$, and consider the cube
\[
{\scriptsize \xymatrix@R=.66em@C=.66em{
& Y_{-\infty} \ar[rr]\ar[dd]\ar[dl]&& Y_0\ar[dd]\ar[dl]\\
 Y_1 \ar[rr]\ar@{=}[dd]&&  Y_2\ar@{=}[dd]\\
& Y_1 \ar[rr]&& Y_2\\
 Y_1 \ar[rr]&&  Y_2\ar@{=}[ur]
}}
\] 
it is easy to see that its {\tt U} = {\tt B}, {\tt F} and {\tt D} faces are cartesian and that the west arrow of its {\tt U} face and the west arrow of its {\tt L} face are isomorphic. The latter arrow lies in $\sS_L$ since the class is 3-for-2, and we can conclude. \qedhere
\end{enumerate}
\end{proof}
\begin{rmk}
It would be really tempting to prove that the three conditions above are in turn equivalent to the following statement, as in categories:
\begin{quote}
\textsc{cf}\oldstylenums{4}) The sub\hyp{}prederivator $\sS_L\subseteq \sD^\bbtwo$ of $L$-locals (\autoref{elle-locals} adapted to prederivators) admits a right calculus of fractions.
\end{quote}
It is in fact quite easy to show that any of the three conditions implies \textsc{cf}\oldstylenums{4}. The equivalence is of course true in case $\D$ is a discrete prederivator; it seems to be a difficult task to prove that condition \textsc{cf}\oldstylenums{4} is sufficient to entail the left exactness of $L$ (and to build a concrete presentation of the prederivator $\D\llbracket \sS_L^{-1}\rrbracket$).
\end{rmk}
\begin{lem}\label{lem:cons-refl}
Let $F\colon\D\to\E$ be a conservative morphism of derivators and let $u\colon A\to B$ be fully faithful. If $F$ preserves left Kan extensions along $u$, then $F$ reflects left Kan extensions along $u$.
\end{lem}
\begin{proof}
Since $u\colon A\to B$ is fully faithful, the same is true for the left Kan extension functors $u_!$ in $\D$ and $\E$, and we can hence characterize the respective essential images by the invertibility of the counits $\varepsilon\colon u_! u^\ast\to\id$. Given a diagram $X\in\D(B)$ by \cite[3.10]{groth:revisit} there is a commutative diagram
\[
\xymatrix{
u_!Fu^\ast\ar@{=>}[r]^-\cong\ar@{=>}[d]_-{\gamma^{-1}}^\cong&Fu_!u^\ast\ar@{=>}[d]^-{F\varepsilon}\\
u_!u^\ast F\ar@{=>}[r]_-{\varepsilon F}&F.
}
\]
in which the unlabeled morphism is the canonical mate expressing the compatibility of $F$ with $u_!$. Since $F$ is assumed to preserve left Kan extensions along $u$, this canonical mate is an isomorphism. Thus, $FX$ lies in the essential image of $u_!\colon\E(A)\to\E(B)$ if and only if $\varepsilon F$ is an isomorphism if and only if $F\varepsilon$ is an isomorphism. Since $F$ is conservative this is the case if and only if $\varepsilon\colon u_! u^\ast X\to X$ is an isomorphism which is to say that $X$ lies in the essential image of $u_!\colon\D(A)\to\D(B)$.
\end{proof}

\begin{cor}\label{cor:cons-refl}
Let $F\colon\D\to\E$ be a conservative morphism of derivators and let $A\in\cCat$. If $F$ preserves colimits of shape $A$, then $F$ reflects colimits of shape $A$.
\end{cor}
\begin{proof}
Since $F$ preserves colimits of shape $A$ if and only if it preserves left Kan extensions along the fully faithful functor $i_A\colon A\to A^\rhd$ (\cite[Prop.~3.9]{groth:revisit}), this statement is immediate from \autoref{lem:cons-refl}.
\end{proof}

\begin{cor}\label{cor:adj-refl-colim}
\begin{enumerate}
\item \label{cor:adj-refl-colimI} Co-reflections of derivators preserve and reflect colimits.
\item \label{cor:adj-refl-colimII} Reflections of derivators preserve and reflect limits.
\item \label{cor:adj-refl-colimIII} Equivalences of derivators preserve and reflect limits and colimits.
\end{enumerate}
\end{cor}
\begin{proof}
This is immediate from \autoref{cor:cons-refl}.
\end{proof}
\begin{rmk}
Let $L\colon\D\rightleftarrows\E\colon R$ be a reflection of derivators.
\begin{enumerate}
\item The left adjoint $L$ preserves left Kan extensions.
\item The left adjoint $L$ preserves right Kan extensions of diagrams in the essential image of $R$.
\item The right adjoint $R$ preserves and reflects right Kan extensions.
\end{enumerate}
It only remains to verify the second claim, and for that purpose let $u\colon A\to B$ be in $\cCat$. In order to show that the canonical mate $Lu_\ast R\to u_\ast LR$ is invertible, it suffices to consider the following diagram
\[
\xymatrix{
LRu_\ast\ar[r]^-\sim\ar[d]_-\varepsilon^-\sim&Lu_\ast R\ar[r]&u_\ast LR\ar[d]^-\varepsilon_-\sim\\
u_\ast\ar[rr]_-\id&&u_\ast
}
\]
which commutes by functoriality of canonical mates and \cite[3.11]{groth:revisit}.
\end{rmk}
\begin{defn}\label{creates}
Let $F : \D \to \E$ be a morphism of prederivators; we say that $F$ \textbf{creates homotopy left} (\textbf{resp., right}) \textbf{Kan extensions along} $u : A\to B$ if for every $X\in \D(B)$, 
\begin{itemize}
	\item the object $X$ lies in the essential image of $u_!^\D$ (resp, $u_*^\D$) if and only if $F_BX$ lies in the essential image of $u_!^\E$ (resp, $u_*^\D$);
	\item $F$ preserves homotopy left (resp., right) Kan extensions along $u$.
\end{itemize}
\end{defn}

\subsection{Orthogonality and co/reflections}
The classical theory motivates a deeper glance at the interaction between reflection of (pre)derivators and the orthogonality relation that can be defined for the objects of $\PDer$; such theory has been used in \cite{tderiv} to sudy \emph{$t$-structures} on stable derivators, finding a derivator analogue of the main theorem in \cite{tstructures}, that ``$t$-structures are normal torsion theories''.

We first recall the definition of \emph{coherent orthogonality} from \cite{factonder}: let $\sD$ be a prederivator satisfying axiom \ref{derax:3}. Let us consider the string of adjoints
\[
\xymatrix{
	\D(\bbtwo) \ar@<-9pt>[r]_{\pt_*}\ar@<9pt>[r]^{\pt_!} & \D(e)\ar[l]|{\pt^*} 
}
\]
where $\pt_!$ (resp. $\pt_*$) sends a coherent diagram $X=\left[ \var{X_0}{X_1} \right]\in\D(\bbtwo)$ into its target object $X_1$ (resp. source object $X_0$) in the base of the derivator.
Given this, we define coherent orthogonality as follows.
\begin{defn}[coherent orthogonality]\label{cohorthogo}
Two coherent diagrams $X,Y\in\D(\bbtwo)$ are called (\textbf{coherently}) \textbf{orthogonal} if the unit morphism $\eta_X  :  X \to \pt^*\pt_!X$ becomes invertible once the functor $\D(\bbtwo)(\firstblank,Y)$ is applied, \ie if the arrow
\[
\D(X_1,Y_0) = \D(\bbtwo)(\pt_!X, \pt_*Y)\cong \D(\bbtwo)(\pt^*\pt_!X,Y) \xto{\D(\bbtwo)(\eta_X,Y)} \D(\bbtwo)(X,Y)
\]
is an isomorphism.
\end{defn}
The orthogonality relation is denoted $X \corth Y$ and defines the condition that ``every commutative square having $X$ on the left and $Y$ on the right admits a unique filler, and coherently so'' in the context of derivators.
This paves the way to the following definition, based on the fact that, dealing with classical orthogonality, we can blur the distinction between objects and their initial or terminal arrows. 
\begin{defn}[co/locality and orthogonality]
Let $X,Y\in\D(\bbtwo)$, $X_1, Y_0 \in\D(e)$.
\begin{enumerate}
	\item We say that $Y_0$ is $X$\textbf{-local}, or $X\corth Y_0$, if $X\corth \left[ \var{Y_0}{1} \right]$;
	\item We say that $X_1$ is $Y$\textbf{-colocal}, or $X_1 \corth Y$, if $\left[ \var{\emptyset}{X_1} \right]\corth Y$;
	\item We say that $X_1$ and $Y_0$ are mutually orthogonal, and write $X_1\corth Y_0$, if $X_1$ is $\left[ \var{Y_0}{1} \right]$-colocal, or equivalently $Y_0$ is $\left[ \var{\emptyset}{X_1} \right]$-local.
\end{enumerate}
This last condition means in particular that $\D(e)(X_1, Y_0)$ is reduced to a singleton.
\end{defn}
Notice that, in the above notation, $X_1\corth Y_0$ if and only if $\D(e)(X_1,Y_0)=0$ and so, depending on the context, we sometimes also use the more common notation $X_1\perp Y_0$ to mean the same as $X_1\corth Y_0$. Certain slight abuses of notation are now straightforward to understand: we can define orthogonality, as well as co/locality, with respect to a chosen \emph{class} of $\{X_\alpha\}_{\alpha\in A} \in \D(\bbtwo)$ and this gives the usual Galois connection
\[
\lcorth{(\firstblank)} \dashv (\firstblank)^{\corth},
\]
which also allows to speak about the pairs $(\sS, \sS^{\corth})$ and $(\lorth{\sS},\sS)$ \emph{generated} by a sub-prederivator $\sS$. The notion of coherent orthogonality is used in \cite{factonder} to lay the foundation of the theory of factorization systems on derivators, and then a theory of coherent $t$-structures as a consequence (see \cite{Fiorenza2016} for the characterization of $t$-structures as ``normal torsion theories''). A general survey of the main features of factorization systems will be the subject of a subsequent work; for the moment we only record that \cite{fs-on-der} already contains a rather general result, since it characterizes certain factorization systems as algebra structures for the ``squaring'' 2-monad of \cite{Korostenski199357}.
\begin{defn}[pre/factorizations and crumblings]
A \textbf{derivator prefactorization system} (\dpfs for short) on a derivator $\D$ is defined to be a pair $\fF=(\sE, \sM)$ such that $\sE = \lcorth{\sM}$ and $\sM = \sE^{\corth}$. A \textbf{derivator factorization system} (\dfs for short) on a derivator $\D$ is defined to be a pair $(\fF,\Psi_\fF)$ where $\fF$ is a \dpfs and $\Psi_\fF$ is a \emph{functorial factorization morphism}, namely (\cite[Def. 3.16]{factonder}) an equivalence of derivators $\Psi : \D_{\mathfrak F}\to \D^\bbtwo$ having domain those $X\in \D^J([2])$ such that $X_{(0,1)} \in\sE(J)$ and $X_{(1,2)}\in\sM(J)$. More generally, if $\sS$ is a sub-prederivator of $\D$, we call $\sS$-\textbf{crumbling} a prefactorization system $(\sE,\sM)$ with a functorial factorization morphism which is an equivalence restricted to $\sS$; of particular importance for us is the case when $\sS$ is the essential image of $\pt_*$; in that case, somewhat sloppily, we say that ``there are factorizations of all terminal arrows'' and we speak about $\tau$-crumbling factorization.
\end{defn}
With \autoref{cohorthogo} in hand, we can prove the derivator analogue of \cite[5.4.4]{Bor1}:
\begin{prop}\label{refl_are_ortho_classes}
Let $L\dashv R  : \D \rightleftarrows \sE$ be a reflection of the derivator $\D$. Consider the prederivator $\sS_L\subseteq \D^\bbtwo$, where $f\in \sS_L(J)$ if and only if $L(f)$ is an isomorphism. Then the following conditions are equivalent for a given $X\in\D^J(\bbtwo)$:
\begin{enumerate}
	\item $X$ lies in the essential image of $R^J_\bbtwo  : \sE^J(\bbtwo) \to \D^J(\bbtwo)$;
	\item $f\corth X$ for every $f\in \sS_L^J(e)$;
	\item given $Y\in\D(J)$, consider $\widetilde{\xi}_{J,Y}\in \D^J(\bbtwo)$ such that $\dia_\bbtwo(\widetilde{\xi}_{J,Y}) : Y\to R_JL_JY$ is the unit of our reflection, then $\widetilde\xi_{J,Y}\corth X$.
\end{enumerate}
Moreover, we can intrinsically characterize the sub\hyp{}prederivator of $L$-locals (\autoref{elle-locals} adapted to prederivators); the following conditions are equivalent for a diagram $X \in \D^J(\bbtwo)$:
\begin{enumerate}
	\item[(1')] $X\in\sS_L^J(e)$;
	\item[(2')] $X\corth Y_0$ for every $Y_0\in \sE^J(e)$;
	\item[(3')] $X \corth Y$ for every $Y\in\sE^J(\bbtwo)$.
\end{enumerate}
\end{prop}
\begin{proof}
A coherent morphism $X\in \D(\bbtwo)$ is an isomorphism (i.e., by definition, $\dia_{\bbtwo}(X)$ is an isomorphism in $\D(e)$) if and only if $\eta_X : X\to \pt^*\pt_!X$ is an isomorphism. In fact, $\eta_X$ is an isomorphism in $\D(\bbtwo)$ if and only if both $(\eta_X)_0$ and $(\eta_X)_1$ are isomorphisms in $\D(e)$; now, $(\eta_X)_1$ is conjugated to the identity of $X_1$, so it is always an isomorphism, while $(\eta_X)_0$ is conjugated to $\dia_\bbtwo X$, so it is an isomorphism exactly when $X$ is so. With this characterization, let us verify the equivalence of statements (1--3):

\smallskip
(1)$\Rightarrow$(2). Let $Y\in \sE^J(\bbtwo)$ be such that $X\cong R^J_\bbtwo Y$. Given $f\in \sS_L^J(\bbtwo)$, we should verify that the canonical map 
\[
\D^J(\bbtwo)(\pt^*\pt_!f,R^J_\bbtwo Y)\cong\D^J(\bbtwo)(\pt^*\pt_!f,X)\to \D^J(\bbtwo)(f,X)\cong \D^J(\bbtwo)(f,R^J_\bbtwo Y).
\]
is an isomorphism.  By adjointness, we can equivalently verify that the following map is an isomorphism: 
\[
\E^J(\bbtwo)(\pt^*\pt_!L^J_\bbtwo f, Y)\cong \E^J(\bbtwo)(L^J_\bbtwo\pt^*\pt_!f, Y)\to \E^J(\bbtwo)(L^J_\bbtwo f, Y).
\]
where the first isomorphism holds since $L^J_\bbtwo$ is a morphism of derivators (so it commutes with $\pt^*$) and it is a left adjoint (so it commutes with $\pt_!$). Now, the above map is induced by the canonical map $\eta_{L^J_\bbtwo f} : L^J_\bbtwo f\to \pt^*\pt_!L^J_\bbtwo f$, which is an isomorphism if and only if $L^J_\bbtwo f$ is an isomorphism, but this is true by hypothesis. 

\smallskip
(2)$\Rightarrow$(3). It is clear that $L_J$ sends the unit $Y\to R_JL_JY$ to an isomorphism, that is, $L_J(\dia_\bbtwo(\widetilde{\xi}_{J,Y}))$ is an isomorphism. Since in a derivator isomorphisms can be checked pointwise, this means that $L_J(\widetilde{\xi}_{J,Y})$ is an isomorphism, that is, $\widetilde{\xi}_{J,Y}\in \sS_L^J(e)$, so the thesis follows.

\smallskip
(3)$\Rightarrow$(1). It is enough, assuming (3), to show that the unit $\xi^J_{\bbtwo,X} : X\to R_{\bbtwo}^JL_{\bbtwo}^JX$ is an isomorphism or, equivalently, that $(\xi^{J}_{\bbtwo,X})_0=\xi_{J,X_0} : X_0\to R_JL_JX_0$ and $(\xi^{J}_{\bbtwo,X})_1=\xi_{J,X_1} : X_1\to R_JL_JX_1$ are both isomorphisms.  Let $i=0,1$, and let $\widetilde{\xi_{J,X_i}}$ be an object lifting $\xi_{J,X_i}$, so by (3), $\widetilde{\xi_{J,X_i}}\corth X$. is clearly an isomorphism 
\end{proof}
We now draw a definition parallel to \autoref{choric_ref}, and in \autoref{refl_are_fs} we prove the equivalence between the two notions. This is the best approximation to the classical result \cite{CHK} connecting reflective subcategories $\cA \subseteq \cC$ with reflective prefactorization systems on $\cC$.
\begin{defn}[Choric factorization]\label{choric_fs}
Let $\D$ be a prederivator, $\fF=(\sE,\sM)$ a pair of sub pre-derivators of $\D^{\bbtwo}$ and let $\D_\fF\colon \Dia\op\to \Cat$ be a pre-derivator such that $\D_{\fF}(I)\subseteq \D(I\times [3])$ is the full subcategory spanned by those $X\in  \D(I\times [3])$ such that $X_{(0,1)}\in \sE(I)$ and $X_{(1,2)}\in \sM(I)$. Denote by 
\[
\label{composition_morphism}
\Psi_\fF\colon \D_\fF\longrightarrow \D^{\bbtwo}
\]
the restriction of the morphism of derivators $\D^{[3]}\to \D^{\bbtwo}$ induced by $(0,2)$.

In this notation, $\fF$ is said to be \emph{choric} if $\Psi_{\fF}$ is fully faithful. More generally, if $\sS$ is a sub-prederivator of $\D$ we say that $\fF$ is $\sS$-choric if $\Psi_{\fF}$ is fully faithful when restricted to $\sS$; particularly important for us is the case when $\sS$ is the essential image of $\pt_*$; in this case we speak about $\tau$-choric \dpfs.
\end{defn}
As we observe in \cite{tderiv}, the \dpfs induced by a pair $(\mathbb E_F,\mathbb M_F)$ is a choric \dfs, and any choric \dfs on $\D$ arises this way from an Eilenberg\hyp{}Moore factorization. So, we are currently unable to exhibit an example example of non-choric \dpfs.
\begin{rmk}
Remind the notation of \autoref{choric_ref}. In the representable case the functor
\[
\D \longleftarrow \chor_L  = \{ X\in\D^{\bbtwo} \mid X_1\in \sM/*\} \subseteq \D^{\bbtwo}
\]
always admits a left adjoint, given by a choice, for a coherent morphism $X\in \cC^{\bbtwo}$, of a dotted arrow in
\[
\xymatrix{X_0\ar[r]^X\ar[d] & X_1\ar[d] \\
RX_0 \ar@{.>}[r] & RX_1.}
\]
Strict orthogonality, and the uniqueness of such a dotted arrow, entail that this choice is in fact unique. But in general, even when liftable to a coherent diagram in $\D^{\bbtwo}(\square)$, the incoherent diagram in $\D(\bbtwo)^\square$ can't be lifted \emph{uniquely}, and additional conditions on a family of functors $\D(I)\to \chor_L(I)$ that exist separately must be imposed in order to ensure that these are the components of a pseudonatural transformation $\D\to \chor_L$.%: in a choric \dpfs these morphisms ``sing together'' giving the desired adjoint.
\end{rmk}
\begin{prop}\label{refl_are_fs}
Given a derivator $\D$, there exists an equivalence between 
\begin{itemize}
	\item the posets of choric reflections of $\D$ (\autoref{choric_ref}) and
	\item the poset of reflective, $\tau$-choric \dpfs on $\D$.
\end{itemize}
\end{prop}
\begin{proof}
Let $L\dashv R$ be a choric reflection, where $R : \E \hookrightarrow \D$; we define $\sS_L$ to be the sub-prederivator of $L$-locals, as we did elsewhere. Our aim is to prove that the \dpfs right-generated by $\sS_L$ is a $\tau$-choric, reflective \dpfs on $\D$. It is obvious that $(\sS_L, \sS_L^{\corth})$ is reflective, as $\sS_L$ is a 3-for-2 class. It is also obvious that is it is a \dpfs, so we are only left to prove that it is $\tau$-choric.

Now, every terminal arrow $\left[\var{X}{*}\right]$ can be $(\sS_L, \sS_L^{\corth})$-factored: if the reflection is choric, there is a well-defined way to attach to $X\in\D$ its reflection $\left[\var{X}{LX}\right]$, and the terminal arrow $\left[\var{LX}{*}\right]$ lies in $\sS_L^{\corth}$; now, this latter arrow lies in $\E(e) \subseteq (\lorth{\E(e)^\bbtwo})^\perp$, so we have the result, and the unit arrow is evidently inverted by $L$.

Conversely, let $(\E,\M)$ be a $\tau$-choric \dpfs on $\D$. We define $\M/*$ to be the class of objects $X\in\D(e)$ such that $t_*X = \smat{X\\ \downarrow\\{*}}$ lies in $\M$. We have to prove that sending $X\in\D(J)$ into $(\Psi_\fF(0_*X))_1$ determines a reflection of $\D$, or more precisely that the composition
\[
L_\fF:\D \xto{0_*} \D^\bbtwo \xto{\Psi_\fF^{-1}} \D_\fF \xto{1^*} \D
\]
when corestricted to its essential image, works as left adjoint to the inclusion $\M/*\hookrightarrow \D$ (and so in particular this essential image coincides with $\M/*$). But this is straightforward, as $L_\fF$ acts on incoherent diagrams as follows:
\[
X \mapsto
\left[ \var{X}{*} \right]
\mapsto
\smat{X &&\\
&\searrow&\\
&&L_\fF X\\
&\swarrow& \\
{*}&&}
\mapsto L_\fF X
\]
% \todo[inline]{PROOF MISSING}
so that the object $L_\fF X$ falls onto $\M/*$.
\end{proof}
\begin{prop}\label{itsafact}
Let $L\dashv R$ be a choric reflection of derivators; its associated \dpfs $(\sS_L, \sS^{\corth}_L)$ is then $\tau$-choric, and if $\sS_L$ is pullback stable then it is also a \dfs.
\end{prop}
\begin{proof}
We prove that $(\sE,\sM)=(\sS_L, \sS_L^{\corth})$ is a factorization. Let $i : \DrawCube{F} \to \DrawCube{DRF}$ and $j : \DrawCube{DRF}\to \DrawCube{BLDRUF}$ be the canonical inclusions of the {\tt F} face and of the {\tt DRF} faces. Given $X\in \D(\bbtwo)$ we consider the coherent morphism $L^\sharp X \to j_*i_!L^\sharp X$ whose incoherent underlying diagram looks like
\[
{\scriptsize
\xymatrix@R=.66em@C=.66em{
& P \ar[rr]\ar[dd]&& RL X_0\ar[dd]\\
 X_0 \ar[ur]\ar[rr]\ar[dd]&&  RL X_0\ar[dd]\ar@{=}[ur]\\
& X_1 \ar[rr]&& RL X_1\\
 X_1 \ar@{=}[ur]\ar[rr]&&  RL X_1\ar@{=}[ur]
}}
\]
Now $RLX\in\D(\bbtwo)$ lies in $\E(e)^{\bbtwo} = (\lorth{\E(e)^{\bbtwo})^\perp} = \sM$, so (since this class is always closed under pullback) also $P\to X_1$ lies in $\sM$. In order t show that the {\tt R} face of the cube is the incoherent diagram associated to the $(\sE,\sM)$-factorization of $X$, we must show that $X_0\to P$ lies in $\sS_L$, \ie that $L$ inverts this morphism: but this follows from the fact that $X_i\to RLX_i$ lies in $\sS_L$ for $i=0,1$ and from the assumption that $\sS_L$ is closed under pullback (recall that $\sS_L$ is a 3-for-2 class).
\end{proof}
% \include*{sections/04-preslim}
%!TEX root = ../localization.tex
\section{Monads and their algebras}
\label{sec:monads}
In this section we introduce some useful terminology and notation about monads on a (pre)derivator. The excellent \cite{lagkas} investigates the theory of monads on derivators, with applications to stable derivators, but also lays the foundation of the general theory in an elegant and readable way. We advise the reader to consult this reference, that we follow quite nearly (in particular \autoref{lagkas_def} and the previous discussion on the derivators of algebras for a monad comes from there).
\begin{defn}[Monad on a prederivator]
Let $\sD$ be a prederivator. We define a \textbf{monad} on $\sD$ to be a morphism $T \colon \sD \to \sD$ equipped with two natural transformations
\begin{enumerate}
	\item $\mu \colon T \circ T \To T$ (the \textbf{multiplication} of the monad);
	\item $\eta \colon \id \To T$ (the \textbf{unit} of the monad),
\end{enumerate}
satisfying the usual associativity and unitality conditions
expressed by the commutativity of the following diagrams of 2-cells: the compatibility of $\mu,\eta$ with the structure of $T$ as a pseudonatural transformation, as well as the associativity and unitality constraints, will usually remain hidden; %for the record, we outline that this spells out as the commutativity %\footnote{Non \'e proprio una commutativit\'a... stiamo dicendo che due due 2-celle sono uguali; \textcolor{red}{In realtà è proprio una richiesta di commutatività: un diagramma strettamente commutativo è riempito da una 2-cella identica; \emph{questo} diagramma di 2-celle è riempito da una 3-cella identica.}} 
%of the following diagrams, for any $u\colon I\to J$ in $\Cat$:
% \begin{gather*}
% \xymatrix@R=.45cm@C=.45cm{
% & \D(J) \ar@{}[d]|{\Downarrow\mu}\ar[dr]^{T_J}&&&&& \D(J) \ar[dd]|{u^*}\ar[dr]^{T_J}\\
% \D(J) \ar@{}[ddrr]|{\Swarrow\tau_u}\ar[dd]^{u^*}\ar[rr]_{T_J}\ar[ur]^{T_J}&& \D(J)\ar[dd]^{u^*}\ar@{}[ddrrr]|{=} &&& \D(J)\ar[dd]_{u^*} \ar[ur]^{T_J}\ar@{}[dr]|{\Swarrow\tau_u}&& \D(J)\ar[dd]^{u^*}\ar@{}[dl]|{\Swarrow\tau_u}\\
% &&&&&&\D(I)\ar[dr]^{T_I}\ar@{}[d]|{\Downarrow\mu}\\
% \D(I) \ar[rr]_{T_I}&& \D(I) &&& \D(I)\ar[ur]^{T_I}\ar[rr]_{T_I} && \D(I)
% }\\
% \xymatrix@R=1.5cm@C=1.5cm{
% \D(J) \rtwocell^{\id}_{T_J}{\eta} \ar[d]_{u^*} & \D(J)\ar[d]^{u^*}\ar@{}[dl]|(.6){\Swarrow\tau_u}\ar@{}[dr]|=  & \D(J) \ar[r]^{\id}\ar[d]_{u^*}& \D(J)\ar[d]^{u^*} \\
% \D(I) \ar[r]_{T_I} & \D(I) & \D(I) \rtwocell^{\id}_{T_I}{\eta} & \D(I)
% }
% \end{gather*}
the relevant diagram can be easily drawn and translated in equational form.
\end{defn}
\begin{rmk}\label{monad_on_each}
It's easy to see that a monad on $\sD$ induces a monad on each category $\sD(J)$, whose multiplication and unit are the components of the modifications $\mu,\eta$ respectively. The fact that the assignment $J \mapsto \sD(J)$ lifts to a suitable category of \emph{categories with monad} $\text{Mon}_l(\Cat)$ (this terminology is better than that in \cite{lagkas}, where these are called \emph{monadic} categories) is a coherence request that can be packed in the following lifting criterion (see \cite{lagkas} again):
\[
\xymatrix{
& \text{Mon}_l(\Cat)\ar[d]^{\boldsymbol U}\\
	\Dia\op \ar[r]_{\D}\ar@{.>}[ur]^{\tilde{\D}} & \Cat
}
\]
(it is possible to obtain $\text{Mon}_l(\Cat)$ as a suitable Grothendieck 2-construction, see \cite{groth2011monoidal}.
\end{rmk}
Of course, a similar strategy yields the definition of an \emph{algebra} for the monad~$T$: \autoref{monad_on_each} above entails that we can consider, for every object $J\in\Cat$, the category of \emph{algebras} for the monad $T_J$ on $\D(J)$. The next result states that all these categories glue together to form a derivator which is the derivator of $\bsT$-algebras on $\D$.
\begin{defn}[\textsc{em}-algebra for a monad]\label{lagkas_def}
Let $\bsT$ be a monad on the prederivator $\D$; the assignment $J\mapsto \D(J)^{T_J}$ defines a prederivator $\D^{\bsT}$ which is the \textbf{prederivator of $\bsT$-algebras} for $\D$, or the \textbf{Eilenberg-Moore prederivators}. Each of the free-forgetful adjunctions $F^{\bsT}_J\dashv U^{\bsT}_J$ glue together as components of an adjunction of prederivators $F^{\bsT}\dashv U^{\bsT}$.
\end{defn}
We now define a 2-category whose 0-cells are the monads over $\D\in\PDer$ (this is of course a particular instance of the 2-category of monads in a 2-category $\cK$, whose 0-cells are the monads over $K\in\cK$).
\begin{defn}[The 2-category of monads in $\PDer$]
Given pairs $(\sE,T)$, $(\sD,S)$ where $\sE,\sD$ are prederivators endowed with monads $T,S$, a \emph{morphism of monads} $(F,\sigma) : (\sE,T)\to (\sD,S)$ is a pair where $F\colon \sE \to \sD$ is a morphism of prederivators, and $\sigma \colon SF \To FT$ is a 2-cell that fills the square
\[
\xymatrix{
	\sE \ar[r]^F\ar[d]_T& \sD \ar[d]^S \ar@{}[dl]|{\Swarrow\sigma}\\
	\sE \ar[r]_F & \sD.
}
\]
This pair is such that the following diagrams of 2-cells commute in $\PDer$:
\[
\begin{tikzcd}
SSF \arrow[r, "\mu_S * F"] \arrow[d, "\sigma\boxminus \sigma"'] & SF \arrow[d, "\sigma "] & F \arrow[r, "F*\eta_T"] \arrow[d, "\eta_S*F"'] & FT \\
FTT \arrow[r, "F * \mu_T"'] & FT & SF. \arrow[ru, "\sigma"] & 
\end{tikzcd}
\]
A 2-cell $\theta : (F,\sigma) \Tto (G,\tau)$ consists of a 2-cell $\theta : F\To G$ in $\PDer$ such that $\tau\circ (S*\theta)= (\theta *T)\circ \sigma$. 

This defines the \emph{2-category of monads in} $\PDer$ that we denote $\Mnd(\PDer)$.
\end{defn}
\begin{rmk}
As observed in \cite{lagkas}, it is easy to see that if $\bsT$ is a monad on $\D$ then each $u^* : \D(J)\to \D(I)$ induced from $u : I\to J$ becomes a strong monad morphism (\ie a monad morphism where $\sigma$ is invertible; the $\sigma$ here is of course $\gamma_{T,u}$ of \autoref{the_coherence}).
\end{rmk}
\begin{rmk}
Note that there exists a 2-functor $\Mnd(\PDer)\to \PDer$ that sends the monad $\bsT$ to its domain and projects $(F,\sigma)$ to $F$. We denote, with a slight imprecise notation, the fiber over $\D$ as $\Mnd(\D)$, the sub-2-category of monads $T : \D\to \D$ on a fixed domain, monad morphisms and monad 2-cells.
\end{rmk}
The following remark is the content of 2.10 and 2.11 on \cite{lagkas}.
\begin{rmk}\label{lax_limi}
A monad on $\sD$ can equivalently be defined as a lax functor $T : * \to \PDer$ from the terminal 2-category. In this picture, as we recall below, the Eilenberg-Moore category of $\bsT$ coincides with the lax limit of $T$.
\end{rmk}
\begin{rmk}
Let $\bsT = (T,\eta,\mu)$ be a monad on a derivator $\sE$; let $\text{Splt}(\bsT)$ the category of those adjoint pairs $F\dashv G$ such that $GF=T, G\varepsilon F =\mu$. The terminal object in this category is called the  ``free-forgetful adjunction''; $\text{Splt}(\bsT)$ can be regarded as the category of cones for $\bsT$, when it is regarded as a lax functor, and the Eilenberg-Moore category of $\bsT$ as its lax limit.
\end{rmk}
\begin{prop}
The 2-category $\PDer$ of prederivators admits the construction of algebras in the sense of \cite{Street1972}. %Moreover, the object of algebras for $\sD\in\PDer$ can be characterized as the lax limit of the monad, regarded as a lax functor $T \colon \bbone \to \PDer$.
\end{prop}
\begin{proof}
Recall that a 2-category ``admits the construction of algebras'' if the canonical functor $\PDer \to \Mnd(\PDer)$ has a right adjoint; unraveling the definition of a monad morphism, we see that a morphism $(F,\sigma) \colon (\sE,\id_\sE)\to (\sD,T)$ in $\Mnd(\PDer)$ is precisely a morphism $\sE \to \D^{\bsT}$. %It remains to check that this defines a 2-adjoint: this is done in \cite[???]{lagkas}.
\end{proof}
% \begin{rmk}
% A more conceptual, but non-constructive strategy would have been to check that $\PDer$ satisfies the sufficient conditions for a 2-category to admit the construction of algebras given in \cite{kelly1974review}: the strategy here is to show that
% 	\begin{itemize}
% 		\item $\PDer$ has finite limits when regarded as a 1-category;
% 		\item the functor $\PDer(\sD,\firstblank) \colon \PDer \to \Cat$ commutes with these finite limits;
% 		\item $\PDer$ has cotensors with all finite categories. 
% 	\end{itemize}
% 	The last condition comes for free from the fact that the shifting operation gives $\Dia$-cotensors to $\PDer$, so we're already set for a big enough choice of $\Dia$. The first condition comes from the fact that $\PDer = [\Dia\op,\Cat]$ inherits all co/limits existing in $\Cat$. The second point is an easy computation with coends, remembering that the hom-category $\PDer(\sD,\sE)$ is the pseudo-end \cite{bozapalides1975fins,bozapalides1977finsgen}
% 	\[
% 		\oint_J \Cat(\sD(J), \sE(J)).
% 	\]
% \end{rmk}
% This is of course a rather terse characterization for the existence of the construction of algebras, that relies on Street's formal theory of monads \cite{Street1972}. 

In this discussion we are guided by the principle that as cumbersome as it may seem, the internal category theory of prederivators cannot stray much from the category theory of $\CAT$, given the tight relation between the two objects:\footnote{A prederivator can be intuitively represented as a ``variable category'' in much the same way in which a presheaf can be thought as a variable sets; this intuition is taken further in \cite{street1981conspectus}, in the case of \emph{pseudo}functors $\cC\to \CAT$ where $\cC$ is a small bicategory.} this suggests that there may exist different way to re-enact monad theory inside $\PDer$, more reminiscent of 1-category theory. Fortunately, it turns out that this is true, and that we can fruitfully borrow many ideas from the treatment of monads in enriched category theory.
\begin{prop}\label{t_modules}
Let $\bsT$ be a monad on the prederivator $\sA$, and $S : \sC \to \sA$ a morphism of prederivators; then $S$ admits a lifting to a morphism of prederivators $\bar S : \sC \to \sA^{\bsT}$ along the forgetful functor $U^{\bsT} : \sA^{\bsT}\to \sA$ if and only if $S$ is a \emph{$T$-module}, \ie if and only if there exists a 2-cell $\zeta : TS\To S$ such that the diagrams
\[
\xymatrix{
TTS \ar[r]^{T\zeta}\ar[d]_{\mu S} & TS\ar[d]^\zeta \\
TS \ar[r]_\zeta & S
}
\qquad
\qquad
\xymatrix{
S \ar@{=}[dr] \ar[r]^{\eta S} & TS \ar[d]^\zeta \\
& S
}
\]
commute. We denote the subcategory of $T$-modules, inside the whole category $\hom(TS,S)$ of 2-cells $TS\To S$, as $\textsf{Act}(TS,S)$
\end{prop}
Let $F : \sX \rightleftarrows \sY : G$ be an adjunction in $\PDer$; then $G$ is a $GF$-module with $\sigma = G\varepsilon $. As a consequence, if we denote $T=GF$ the associated monad, we have
\begin{cor}\label{comparison_def}
There exists a unique morphism of prederivators $K = K_{F,G} : \sY \to \sX^{\bsT}$, that moreover has the property $KF=F^{\bsT}$. This morphism is called the \textbf{comparison morphism} between $T$-algebras and $\sY$.
\end{cor} 
\begin{rmk}
The assignment $\fS : \Mnd(\sA) \to \PDer_{\!/\sA} : \bsT\mapsto \left[U^{\bsT} : \var{\sA^{\bsT}}{\sA}\right]$ defines a contravariant functor that realizes a bijection (natural in $T$, with respect to monad morphisms)
\[
\PDer(\sX,\sA)(S,\fS(\bsT)) \cong \textsf{Act}(TS,S)
\]
where in the right hand side we denoted the subcategory of $T$-modules $TS\To S$ defined in \autoref{t_modules}. Moreover, if we track the image of $\textsf{Act}(TS,S)$ under the chain of isomorphisms
\[
\PDer(\sX, \sA)(S, \fS(T))\cong \textsf{Act}(TS,S) \subset \hom(TS,S) \cong \hom(T, \Ran_SS)
\]
we get that this corresponds exactly to the subcategory of monad morphisms $T \to \Ran_SS$ between $T$ and the \emph{codensity monad} $\langle S,S\rangle = (\Ran_SS,u_S, \varepsilon_S)$. Thus we established the equivalence
\[
\PDer(\sX, \sA)(S, \fS(\bsT)) \cong \Mnd(\sA)(\bsT,\langle S,S\rangle),
\] that in high-sounding terms can be called \cite{dubuc1970kan,Lack2010} the ``semantic $\dashv$ structure'' adjunction in $\PDer$ between the 2-category $\PDer_{\!/\sA}^{\leaden}$ of those $S : \sX\to \sA$ admitting a codensity monad.\footnote{Read ``heavy $\PDer_{\!/\sA}$'' for the subcategory of functors $\sX\to \sA$ admitting a codensity monad; $\leaden$ is the alchemical token symbolizing lead (or planet Saturn) \cite[p. 153]{historyalchemy}.}
\end{rmk}
It is interesting to single out the monads whose multiplication is an isomorphism; because of their property these are called \emph{idempotent}, and it turns out that they correspond to reflective subcategories of their domain, under the correspondence $\bsT\mapsto \D^{\bsT}$.
\begin{defn-prop}\label{char_of_ref_using_mnds}
(this is the derivator analogue of \cite[4.2.4]{Bor1}) The following conditions are equivalent
\begin{enumerate}[label=\textsc{im}\oldstylenums{\arabic*}), ref=\textsc{im}\oldstylenums{\arabic*}]
	\item \label{idem_mon:2}the counit of the adjunction $F^{\bsT}\dashv U^{\bsT}$ is an isomorphism;
	\item \label{idem_mon:3}the multiplication $\mu \colon T\circ T \To T$ is an invertible modification;
	\item \label{idem_mon:4}For every $T$-algebra $(A,a)$, the structure map $TA\to A$ is an isomorphism.
\end{enumerate}
If any of these conditions is satisfied, $T$ is called an \emph{idempotent monad}.
\end{defn-prop}
\begin{rmk}
According to \autoref{defn:refl-loc}, condition \ref{idem_mon:2} entails that there is a reflection $\D^{\bsT} \rightleftarrows \sD$; % it suffices then to prove \autoref{char_of_ref_using_mnds} to prove also the analogue of the classical characterization of reflective sub-prederivators of $\sD$ as categories of \textsc{em}-algebras for idempotent monads on $\sD$, stated below in \autoref{bij_with_ref}.
we remark that \cite[4.2.4]{Bor1} actually contains another equivalent condition ($1'$): a monad $T$ is idempotent if the forgetful $U^{\bsT} \colon \D^{\bsT} \to \sD$ is fully faithful. In the category of prederivators, we take condition \ref{idem_mon:2} as a \emph{definition} of $U^{\bsT}$ being fully faithful, so the equivalence $(1')\iff (1)$ is true by definition.
\end{rmk}
\begin{proof}
We show the chain of implications $(2)\To (3)\To (1) \To (2)$ and then conclude.

By definition, an invertible modification is such if and only if it has invertible components: since a $T$-algebra is defined to be the result of gluing together all the $T_J$-algebras on $\sD(J)$ we can think of a $T$-algebra as a $\Cat$-indexed family of morphisms $a_J \in\sD(J)(T_JA_J,A_J)$.\footnote{Even though at the end of this proof we will see that the $a_J$ can be chosen naturally in $J$, we must note that this nice behaviour is a consequence of the idempotency of the monad $T$: in a few words, a classical argument shows how the uniqueness of the $T$-algebra structure on an object forces it to be natural.}

A $T$-algebra now satisfies the first two commutativity conditions of this list, and the unit of the monad relates to the multiplication as depicted in the third diagram:
\[
\xymatrix{
	T_J A_J \ar@{}[dr]|(.3){(\textsc{i})}\ar[r]^{a_J}& A_J \\
	A_J\ar[u]^{\eta_J}\ar@{=}[ur] &
}\quad
\xymatrix{
	T_JT_JA_J \ar@{}[dr]|{(\textsc{ii})}\ar[r]^{T_J a_J}\ar[d]_{\mu_J}& T_J A_J \ar[d]^{a_J}\\
	T_JA_J \ar[r]_{a_J}& A_J
}\quad
\xymatrix{
	T_J A_J\ar[r]^{T_J \eta_J}\ar@{=}[dr] & T_JT_JA_J \ar[d]_{\mu_J}& T_JA_J \ar[l]_{\eta_J T_J}\ar@{=}[dl]\\
	& T_J A_J &
}
\]
Assuming $(2)$, from these relations we get that $\eta_J * T_J = T_J * \eta_J = \mu_J^{-1}$ (since this is true on each component, we can safely write $\eta* T = T * \eta = \mu^{-1}$). But then $T_J a_J$ is invertible and equal to $(T_J* \eta_J)^{-1} = (\mu_J^{1})^{-1} = \mu_J$. Finally, naturality for $\eta_J$ gives that
\[
\xymatrix@C=2cm{
	T_J A_J \ar[r]^{\eta_J * T_J A_J}\ar[d]_{a_J} & T_J T_J A_J \ar[d]^{T_J a_J} \\
	A_J \ar[r]_{\eta_{A_J}} & T_J A_J
}
\]
Since $\eta_{A_J}\circ a_J = T_J a_J \circ \eta_J * T_J = \id_{T_JA_J}$, this gives that $a_J$ is forced to be the inverse (hence unique) for $\eta_{A_J}$.

To show that $(3)\To (1)$, and in fact that $(1)\iff (2)$, it suffices to write the components of the counit $F^{\bsT}U^{\bsT} \To 1$: it is evidently diagram II above.
\end{proof}
This result proves part of the following, equivalent characterizations of localizations of a prederivator:
\begin{thm}\label{bij_with_ref}
\cite[5.5.6]{Bor1} There is a bijection between
\begin{enumerate}[label=\oldstylenums{\arabic*}), ref=\oldstylenums{\arabic*}]
	\item \label{refl_char:uno} $\tau$-choric localizations of $\sD$;
	\item \label{refl_char:due} (categories of algebras for) idempotent monads on $\sD$.
	\item \label{refl_char:tre} reflective, $\tau$-choric \dpfs on $\sD$.
\end{enumerate}
This bijection restricts to a bijection between
\begin{enumerate}[label=\oldstylenums{\arabic*}\textsc{r}), ref=\oldstylenums{\arabic*}\textsc{r}]
	\item \label{refl_charstar:uno} left exact $\tau$-choric reflections of $\sD$;
	\item \label{refl_charstar:due} (categories of algebras for) left exact idempotent monads on $\sD$.
	\item \label{refl_charstar:tre} reflective, $\tau$-choric \dpfs on $\sD$, closed under finite Kan extensions.
\end{enumerate}
\end{thm}
\begin{proof}
The equivalence $1\iff 2$ follows from \autoref{char_of_ref_using_mnds}, and the equivalence $2\iff 3$ is the content of \autoref{refl_are_fs}.
% \todo[inline]{PROOF MISSING}

We now show how this bijection restricts when we consider left exact reflections.\footnote{Recall from \autoref{kinds_of_refle} that a \emph{left exact} localization is a localization whose left adjoint commutes with finite right Kan extensions, and that a \emph{finite} right Kan extension is a functor of the form $u_*$ for $u : I\to J$ a functor between homotopy finite categories --categories whose nerve has a finite number of nondegenerate simplices as a simplicial set.} Let $u : I\to J$ be a functor between homotopy finite categories, $\bsT$ and idempotent monad on $\D\in\PDer$, and $F^{\bsT} : \D^{\bsT}\rightleftarrows \D : U^{\bsT}$ the free-forgetful adjunction of \autoref{lagkas_def}. Since $F^{\bsT}$ and $U^{\bsT}$ and $T = U^{\bsT}F^{\bsT}$ are morphisms of derivators, they come equipped with 2-cells $(\gamma_{F^{\bsT},u})_*, (\gamma_{U^{\bsT},u})_*, (\gamma_{T,u})_*$ as in \autoref{tehgammas} such that
\[
\begin{tikzcd}
\D(I) \arrow[d, "F^{\bsT}_I"'] \arrow[r, "u_*"] & \D(J) \arrow[d, "F^{\bsT}_J"] \arrow[shorten <>=3mm, implies,ld, "(\gamma_{F^{\bsT}})_*"'] & \D(I) \arrow[dd, "T^{\bsT}_I"'] \arrow[r, "u_*"] & \D(J) \arrow[dd, "T^{\bsT}_J"] \arrow[shorten <>=3mm, implies,ldd,"(\gamma_{T})_*"'] & U^{\bsT}F^{\bsT}u_*  \arrow[dd, "{U^{\bsT} (\gamma_{F^{\bsT},u})_*}"',implies] \arrow[rrdd, "{(\gamma_{T,u})_*}",implies] &  &  \\
\D^{\bsT}(I) \arrow[d, "U^{\bsT}_I"'] \arrow[r, "u_*"] & \D^{\bsT}(J) \arrow[d, "U^{\bsT}_J"] \arrow[shorten <>=3mm, implies,ld,"(\gamma_{U^{\bsT}})_*"'] &  &  &  &  &  \\
\D(I) \arrow[r, "u_*"'] & \D(J) & \D(I) \arrow[r, "u_*"'] & \D(J) & U^{\bsT}u_*F^{\bsT} \arrow[rr, "{(\gamma_{U^{\bsT},u})_* F^{\bsT}}"',implies] &  & u_*U^{\bsT}F^{\bsT}
\end{tikzcd}
\]
but now the 2-cell $(\gamma_{T,u})_*$ is invertible (\ie $T$ commutes with homotopy finite Kan extensions) if and only if $U^T * (\gamma_{F^{\bsT},u})_*$ is invertible; since $U^{\bsT}$ has each component fully faithful, it is conservative, and this latter condition is true if and only if $(\gamma_{F^{\bsT},u})_*$ is invertible (\ie if $F^{\bsT}$ commutes with homotopy finite Kan extensions).
\end{proof}

\hrulefill

\subsubsection*{Acknowledgements.} The author is supported by the Grant Agency of the Czech Republic under the grant \textsc{P}201/12/\textsc{G}028.

\bibliographystyle{alpha}
\bibliography{allofthem}

\end{document}